\documentclass[letterpaper,11pt]{article}
\usepackage{fullpage,graphicx,psfrag,amsthm,amssymb,subfigure}
\usepackage[margin=1in]{geometry}
\newtheorem{theorem}{Theorem}[section]
\newtheorem{corollary}[theorem]{Corollary}

\newcommand{\cf}{{\it cf.}}
\newcommand{\eg}{{\it e.g.}}
\newcommand{\ie}{{\it i.e.}}

\newcommand{\ones}{\mathbf 1}
\newcommand{\reals}{{\mbox{\bf R}}}
\newcommand{\Diag}{\mathop{\bf Diag}}
\newcommand{\Tr}{\mathop{\bf Tr}}

\title{Fastest Mixing Markov Chain on Graphs with Symmetries}

\author{
Stephen Boyd\thanks{Department of Electrical Engineering,
  Stanford University, Stanford, CA 94305. 
  Email: \texttt{boyd@stanford.edu}.} \and
Persi Diaconis\thanks{Department of Statistics and Department of
  Mathematics, Stanford University, Stanford, CA 94305.} \and
Pablo A. Parrilo\thanks{Department of Electrical Engineering and Computer
  Science, Massachusetts Institute of Technology, Cambridge, MA 02139. 
  Email: \texttt{parrilo@mit.edu}.} \and 
Lin Xiao\thanks{Microsoft Research, 1 Microsoft Way, Redmond, WA 98052. 
  Email: \texttt{lin.xiao@microsoft.com}.}
}

\date{April 27, 2007}

\begin{document}
\maketitle

\begin{abstract}
We show how to exploit symmetries of a graph to efficiently
compute the fastest mixing Markov chain on the graph 
(\ie, find the transition probabilities on the edges to minimize the 
second-largest eigenvalue modulus of the transition probability matrix).
Exploiting symmetry can lead to significant reduction in
both the number of variables and the size of matrices in 
the corresponding semidefinite program, 
thus enable numerical solution of large-scale instances that 
are otherwise computationally infeasible. 
We obtain analytic or semi-analytic results for particular classes
of graphs, such as edge-transitive and distance-transitive graphs.
We describe two general approaches for symmetry exploitation,
based on orbit theory and block-diagonalization, respectively.
We also establish the connection between these two approaches. 

\end{abstract}

\section{Introduction}
In the fastest mixing Markov chain problem,
we choose the transition probabilities on the edges of a graph to
minimize the second-largest eigenvalue modulus of the
transition probability matrix.
In~\cite{BDX:04} we formulated this problem as a convex
optimization problem, in particular as a semidefinite program.
Thus it can be solved, up to any given precision, in polynomial time
by interior-point methods.
In this paper, we show how to exploit symmetries of a graph to make
the computation more efficient.

\subsection{The fastest mixing Markov chain problem}
We consider an undirected graph $\mathcal{G}=(\mathcal{V},\mathcal{E})$ 
with vertex set $\mathcal{V}=\{1,\ldots,n\}$ and edge
set~$\mathcal{E}$ and assume that $\mathcal{G}$ is connected. 
We define a discrete-time Markov chain on the vertices as follows.  
The state at time $t$ will be denoted $X(t) \in \mathcal V$, for
$t=0,1, \ldots $.
Each edge in the graph is associated with a transition
probability with which $X$ makes a transition between the two adjacent
vertices. 
This Markov chain can be described via its transition probability
matrix $P\in\reals^{n\times n}$, where
\[
P_{ij} = \mbox{\bf Prob}\; (\;X(t+1)=j \;|\; X(t)=i\;), \qquad 
i,~j=1, \ldots, n.
\]
Note that $P_{ii}$ is the probability that $X(t)$ stays at
vertex~$i$, and $P_{ij}=0$ for $\{i,j\}\notin\mathcal{E}$ 
(transitions are allowed only between vertices that are linked by
an edge).
We assume that the transition probabilities are symmetric, \ie,
$P=P^T$, where the superscript $T$ denotes the transpose of a matrix.    
Of course this transition probability matrix must also be stochastic:
\[
P\geq 0, \qquad  P\ones=\ones,
\]
where the inequality $P \geq 0$ means elementwise, and $\ones$ denotes
the vector of all ones.

Since $P$ is symmetric and stochastic, 
the uniform distribution $(1/n)\ones^T$ is stationary.
In addition, the eigenvalues of~$P$ are real, and no more than one in
modulus. 
We denote them in non-increasing order
\[
1=\lambda_1(P)\geq\lambda_2(P)\geq\cdots\geq\lambda_n(P)\geq-1.
\]
We denote by $\mu(P)$ the second-largest eigenvalue modulus (SLEM) of
$P$, \ie,
\[
\mu(P) = \max_{i=2,\ldots,n} |\lambda_i(P)| = \max \; \{ \lambda_2(P), \; -\lambda_n(P) \}. 
\]
This quantity is widely used to bound the asymptotic convergence rate
of the distribution of the Markov chain to its stationary
distribution, in the total variation distance or chi-squared distance 
(\eg, \cite{DiS:91,DSC:93}).
In general the smaller $\mu(P)$ is, the faster the Markov chain
converges.
For more background on Markov chains, eigenvalues and rapid mixing, 
see, \eg, the text~\cite{Bre:99}.

In \cite{BDX:04}, we addressed the following problem:
What choice of $P$ minimizes $\mu(P)$?
In other words, what is the fastest mixing (symmetric) Markov chain 
on the graph?
This can be posed as the following optimization problem:
\begin{equation}\label{e-fmmc}
\begin{array}{ll}
\mbox{minimize} &  \mu(P)\\
\mbox{subject to} & 
P\geq 0, \quad  P\ones =\ones, \quad P = P^T\\
&P_{ij} = 0, \quad \{i,j\}\notin \mathcal{E}.
\end{array}
\end{equation}
Here $P$ is the optimization variable, and the graph is the  
problem data.
We call this problem the \emph{fastest mixing Markov chain} (FMMC)
problem. 
This is a convex optimization
problem, in particular, the objective function can be
explicitly written in a convex form 
$\mu(P) = \|P-(1/n)\ones\ones^T\|_2$,
where $\|\cdot\|_2$ denotes the spectral norm of a matrix. 
Moreover, this problem can be readily transformed into a semidefinite
program (SDP):
\begin{equation}\label{e-fmmc-sdp}
\begin{array}{ll}
\mbox{minimize} &  s\\
\mbox{subject to} & 
-sI \preceq P - (1/n)\ones\ones^T \preceq sI\\
&P\geq 0, \quad P\ones =\ones, \quad P = P^T\\
&P_{ij} = 0, \quad \{i,j\}\notin \mathcal{E}.
\end{array}
\end{equation}
Here $I$ denote the identity matrix, and the variables are the 
matrix $P$ and the scalar $s$.  
The symbol $\preceq$ denotes matrix inequality, \ie, $X \preceq Y$
means $Y-X$ is positive semidefinite.

There has been some follow-up work on this problem. 
Boyd, Diaconis, Sun, Xiao (\cite{BDSX:06}) proved analytically that on an $n$-path the fastest mixing chain can be obtained by assigning the same transition probability half at the $n-1$ edges and two loops at the two ends. 
Roch (\cite{Roc:05}) used standard mixing-time analysis techniques (variational characterizations, conductance, canonical paths) to bound the fastest mixing time.
Gade and Overton (\cite{GaO:06}) have considered the fastest mixing problem for a nonreversible Markov chain. Here, the problem is non-convex and much remains to be done.  
Finally, closed form solutions of fastest mixing problems have recently been applid in statistics to give a generalization of the usual spectral analysis of time series for more general discrete data. see \cite{Sal:06}.

\subsection{Exploiting problem structure}
\label{s-structure}
The SDP formulation~(\ref{e-fmmc-sdp}) means that the FMMC problem can
be efficiently solved using standard SDP solvers, at least for small
or medium size problems (with number of edges up to a thousand or so).
General background on convex optimization and SDP can be found in,
\eg, \cite{NeN:94,VaB:96,WSV:00,BTN:01,BoV:04}.
The current SDP solvers (\eg, \cite{Stu:99,SDPT3,SDPA6}) mostly use 
interior-point methods which have polynomial time worst-case
complexity. 

When solving the SDP~(\ref{e-fmmc-sdp}) by interior-point methods, 
in each iteration we need to compute the first and second derivatives
of the logarithmic barrier functions (or potential functions) for 
the matrix inequalities, and assemble and solve a linear system of
equations (the Newton system). 
Let~$n$ be the number of vertices and~$m$ be the number of edges in
the graph (equivalently $m$ is the number of variables in the
problem). 
The Newton system is a set of~$m$ linear equations with~$m$ unknowns.
Without exploiting any structure, the number of flops per iteration in
a typical barrier method is on the order $\max\{m n^3,m^2 n^2,m^3\}$,
where the first two terms come from computing and assembling the
Newton system, and the third term amounts to solving it (see, \eg,
\cite[\S11.8.3]{BoV:04}).
(Other variants of interior-point methods have similar orders of flop
count.)  

Exploiting problem structure can lead to significant improvement of
solution efficiency. 
As for many other problems defined on a graph, sparsity is the most 
obvious structure to consider here.
In fact, many current SDP solvers already exploit sparsity.
However, as a well-known fact, exploiting sparsity alone in
interior-point methods for SDP has limited effectiveness.
The sparsity of $P$, and the sparsity plus rank-one structure of
$P-(1/n)\ones\ones^T$, can be exploited to significantly reduce the
complexity of assembling the Newton system, 
but typically the Newton system itself is dense.
The computational cost per iteration can be reduced to order
$O(m^3)$, dominated by solving the dense linear system 
(see analysis for similar problems in, \eg, \cite{BYZ:00,XiB:04,XBK:07}).

In addition to using interior-point methods for the SDP
formulation~(\ref{e-fmmc-sdp}),  
we can also solve the FMMC problem in the form~(\ref{e-fmmc})
by subgradient-type (first-order) methods.
The subgradients of $\mu(P)$ can be obtained by computing the extreme
eigenvalues and associated eigenvectors of the matrix~$P$.
This can be done very efficiently by iterative methods, specifically
the Lanczos method, for large sparse symmetric matrices (\eg,
\cite{GoV:96,Saa:92}).  
Compared with interior-point methods, subgradient-type methods
can solve much larger problems but only to a moderate accuracy
(they don't have polynomial-time worst-case complexity).
In~\cite{BDX:04}, we used a simple subgradient method to solve the
FMMC problem on graphs with up to a few hundred thousand edges.
More sophisticated first-order methods for solving large-scale
eigenvalue optimization problems and SDPs have been reported in, \eg,  
\cite{HeR:00,BuM:03,Nem:04,LNM:04,Nes:05}.
A successive partial linear programming method was developed in 
\cite{Ove:92}. 

In this paper, we focus on the FMMC problem on graphs with large 
symmetry groups, and show how to exploit symmetries of the graph to
make the computation more efficient. 
A result by Erd{\H o}s and R\'enyi \cite{ErR:63} states that with 
probability one, 
the symmetry group of a (suitably defined) random graph is trivial,
\ie, it contains only the identity element.
Nevertheless, many of the graphs of theoretical and practical
interest, particularly in engineering applications have very
interesting, and sometimes very large, symmetry groups.
Symmetry reduction techniques have been explored in several different
contexts, \eg, dynamical systems and bifurcation theory \cite{GSS:88},
polynomial system solving \cite{Gat:00,Wor:94}, numerical solution of
partial differential equations \cite{FaS:92}, and Lie symmetry
analysis in geometric mechanics \cite{MaR:99}.
In the context of optimization, a class of SDPs with symmetry has been 
defined in \cite{KOMK:01}, where the authors study the invariance 
properties of the search directions of primal-dual interior-point methods.  
In addition, symmetry has been exploited to prune the enumeration 
tree in branch-and-cut algorithms for integer programming \cite{Mar:03}, 
and to reduce matrix size in a spectral radius optimization problem 
\cite{HOY:03}.

Closely related to our approach in this paper, the recent work 
\cite{dPS:07} considers general SDPs that are invariant under the 
action of a permutation group, and developed a technique based on
matrix $*$-representation to reduce problem size.
This technique has been applied to simplify computations in SDP 
relaxations for graph coloring and maximal clique problems 
\cite{DuR:07}, 
and to strengthen SDP bounds for some coding problems \cite{Lau:07}.

For the FMMC problem, we show that exploiting symmetry allows
significant reduction in both number of optimization variables and
size of matrices.  Effectively, they correspond to reducing~$m$
and~$n$, respectively, in the flop counts for interior-point methods
mentioned above.  
The problem can be considerably simplified and is often solvable 
analytically by only exploiting symmetry. 
We present two general approaches for symmetry exploitation,
based on orbit theory \cite{BDPX:05} and block-diagonalization 
\cite{GatermannParrilo}, respectively.
We also establish the connection between these two approaches.

\subsection{Outline}

In~\S\ref{s-analysis}, we explain the concepts of graph automorphisms
and the automorphism group (symmetry group) of a graph.  We show that
the FMMC problem always attains its optimum in the fixed-point subset
of the feasible set under the automorphism group.  This allows us to
only consider a number of distinct transition probabilities that
equals the number of orbits of the edges.  We then give a formulation
of the FMMC problem with reduced number of variables (transition
probabilities), which appears to be very convenient in subsequent
sections.

In~\S\ref{s-analytic}, we give closed-form solutions for
the FMMC problem on some special classes of graphs, namely 
edge-transitive graphs and distance-transitive graphs. 
Along the way we also discuss FMMC on graphs formed by taking
Cartesian products of simple graphs. 

In~\S\ref{s-fmmc-orbit}, we first review the orbit theory for
reversible Markov chains, and give sufficient conditions on
constructing an orbit chain that contain all distinct eigenvalues of
the original chain. 
This orbit chain is usually no longer symmetric but always
reversible. 
We then solve the fastest reversible Markov chain problem on the orbit 
graph, from which we immediately obtain optimal solution to the
original FMMC problem. 

In~\S\ref{s-block-diag}, we review some group representation theory
and show how to block diagonalize the linear matrix inequalities in
the FMMC problem by constructing a symmetry-adapted basis.
The resulting blocks usually have much smaller sizes and repeated
blocked can be discarded in computation. 
Extensive examples in~\S\ref{s-fmmc-orbit} and~\S\ref{s-block-diag}
reveal interesting connections between these two general symmetry
reduction methods.  

In~\S\ref{s-conclusions}, we conclude the paper by pointing out some
possible future work.

\section{Symmetry analysis}
\label{s-analysis}

In this section we explain the basic concepts that are essential in
exploiting graph symmetry, and derive our result on reducing the 
number of optimization variables in the FMMC problem. 

\subsection{Graph automorphisms and classes}
The study of graphs that possess particular kinds of symmetry
properties has a long history.
The basic object of study is the \emph{automorphism group} of a graph,
and different classes can be defined depending on the specific form in
which the group acts on the vertices and edges.

\begin{figure}
\begin{center}
\psfrag{1}[cl]{$1$}
\psfrag{2}[tc]{$2$}
\psfrag{3}[tc]{$3$}
\psfrag{4}[bl]{$4$}
\includegraphics[width=0.55\textwidth]{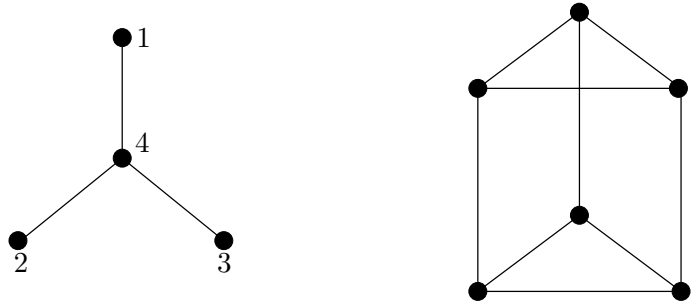}
\caption{The graph on the left side is edge-transitive, but not
  vertex-transitive. The one on the right side is vertex-transitive,
  but not edge-transitive.}
\label{fig:etvt}
\end{center}
\end{figure}

An \emph{automorphism} of a graph
$\mathcal{G}=(\mathcal{V},\mathcal{E})$ is a permutation $\sigma$ of
$\mathcal{V}$ such that $\{i,j\} \in \mathcal{E}$ if and only if
$\{\sigma(i),\sigma(j)\} \in \mathcal{E}$.  The (full)
\emph{automorphism group} of the graph, denoted by
$\mbox{Aut}(\mathcal{G})$, is the set of all such permutations, with
the group operation being composition.  For instance, for the graph on
the left in Figure~\ref{fig:etvt}, the corresponding automorphism
group is generated by all permutations of the vertices $\{1,2,3\}$.
This group, isomorphic to the symmetric group $S_3$, has six elements,
namely the permutations $123\to123$ (the identity), $123\to213$,
$123\to132$, $123\to321$, $123\to231$, and $123\to312$.  Note that
vertex~$4$ cannot be permuted with any other vertex.

Recall that an \emph{action} of a group~$G$ on a set~$\mathcal X$ is a
homomorphism from~$G$ to the set of all permutations of the elements
in~$\mathcal X$ (\ie, the symmetric group of degree $|\mathcal X|$). 
For an element $x\in\mathcal X$, the set of all images $g(x)$, as
$g$ varies through $G$, is called the \emph{orbit} of~$x$. 
Distinct orbits form equivalent classes and they partition the 
set~$\mathcal X$. 
The action is \emph{transitive} if for every pair of elements $x, y \in
\mathcal{X}$, there is a group element $g \in G$ such that $g(x)=y$. 
In other words, the action is transitive if there is only one single 
orbit in~$\mathcal{X}$.

A graph $\mathcal{G} = (\mathcal{V}, \mathcal{E})$ is said to be
\emph{vertex-transitive} if $\mbox{Aut}(\mathcal G)$ acts transitively
on $\mathcal{V}$.
The action of a permutation $\sigma$ on $\mathcal V$ induces an
action on $\mathcal E$ with the rule 
$\sigma(\{i,j\})=\{\sigma(i),\sigma(j)\}$. 
A graph $\mathcal{G}$ is \emph{edge-transitive} if
$\mbox{Aut}(\mathcal G)$ acts transitively on $\mathcal{E}$. 
Graphs can be edge-transitive without being vertex-transitive and 
vice versa; simple examples are shown in Figure~\ref{fig:etvt}. 

A graph is \emph{1-arc-transitive} if given any four vertices
$u,v,x,y$ with $\{u,v\},\{x,y\} \in \mathcal{E}$, there exists an
automorphism $g \in \mbox{Aut}(\mathcal G)$ such that $g(u) = x$ and
$g(v) = y$.
Notice that, as opposed to edge-transitivity, here the ordering of
the vertices is important, even for undirected graphs. 
In fact, a 1-arc transitive graph must be both vertex-transitive and
edge-transitive, and the reverse may not be true. 
The 1-arc-transitive graphs are called \emph{symmetric graphs} in
\cite{Biggs}, but the modern use extends this term to all graphs that
are simultaneously edge- and vertex-transitive.
Finally, let $\delta(u,v)$ denote the distance between two vertices
$u,v \in \mathcal{V}$.
A graph is called \emph{distance-transitive} if, for any four vertices
$u,v,x,y$ with $\delta(u,v) = \delta(x,y)$, there is an automorphism
$g \in \mbox{Aut}(\mathcal G)$ such that $g(u) = x$ and $g(v) = y$.

\begin{figure}
\begin{center}
\setlength{\unitlength}{.02cm}
  \begin{picture}(350,200)(0,30)
    \put(80,120){\framebox(230,40){\textbf{1-arc transitive}}}
    \put(80,200){\framebox(230,40){\textbf{Distance-transitive}}}
    \put(220,40){\framebox(190,40){\textbf{Vertex-transitive}}}
    \put(-30,40){\framebox(190,40){\textbf{Edge-transitive}}}
    \put(185,200){\vector(0,-1){40}}
    \put(220,120){\vector(1,-1){40}}
    \put(150,120){\vector(-1,-1){40}}
\end{picture}
\end{center}
\caption{Classes of symmetric graphs, and their inclusion relationship.}
\label{fig:graphclasses}
\end{figure}
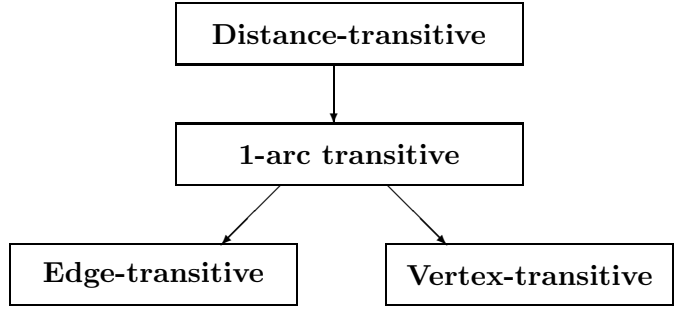

The containment relationship among the four classes of graphs
described above is illustrated in Figure~\ref{fig:graphclasses}. 
Explicit counterexamples are known for each of the non-inclusions.
It is generally believed that distance-transitive graphs have been completely classified. This work has been done by classifying the distance-regular graphs. 
It would take us too far afield to give a complete discussion.
See the survey in \cite[Section 7]{DSC:06}.

The concept of graph automorphism can be naturally extended to
weighted graphs, by requiring that the permutation must also preserve
the weights on the edges (\eg, \cite{BDPX:05}).
This extension allows us to exploit symmetry in more general
reversible Markov chains, where the transition probability matrix is
not necessarily symmetric. 

\subsection{FMMC with symmetry constraints}
\label{s-fmmc-invariant}
A permutation $\sigma\in\mbox{Aut}(\mathcal G)$ can be represented by
a permutation matrix $Q$, where $Q_{ij}=1$ if $i=\sigma(j)$ and
$Q_{ij}=0$ otherwise.
The permutation $\sigma$ induces an action on the transition
probability matrix by $\sigma(P)=Q P Q^T$.
We denote the feasible set of the FMMC problem~(\ref{e-fmmc}) by  
$\mathcal C$, \ie,
\[
\mathcal C = \{ P\in\reals^{n\times n} ~|~ P\geq 0, ~P\ones=\ones,
~P=P^T, P_{ij}=0 ~\mbox{for}~\{i,j\}\notin\mathcal E\}.
\]
This set is invariant under the action of graph automorphism. 
To see this, let $h=\sigma(i)$ and $k=\sigma(j)$.
Then we have 
\[
(\sigma(P))_{hk} = (QPQ^T)_{hk} = \sum_l (QP)_{hl}Q_{kl} = (QP)_{hj}
= \sum_l Q_{hl}P_{lj}=P_{ij}.
\]
Since $\sigma$ is a graph automorphism, we have $\{h,k\}\in\mathcal E$
if and only if $\{i,j\}\in\mathcal E$, so the sparsity pattern of the
probability transition matrix is preserved.
It is straightforward to verify that the conditions $P\geq0$,
$P\ones=\ones$, and $P=P^T$, are also preserved under this action.

Let~$\mathcal F$ denote the fixed-point subset
of~$\mathcal C$ under the action of $\mbox{Aut}(\mathcal G)$; \ie,
\[
\mathcal F = \{ P\in\mathcal C ~|~ \sigma(P)=P,
~\sigma\in\mbox{Aut}(\mathcal G)\}.
\]
We have the following theorem (see also 
\cite[Theorem~3.3]{GatermannParrilo}).
\begin{theorem}\label{t-fmmc-fixed-point}
The FMMC problem always has an optimal solution in
the fixed-point subset $\mathcal F$.
\end{theorem}
\begin{proof}
Let $\mu^\star$ denote the optimal value of the FMMC
problem~(\ref{e-fmmc}), \ie, 
$\mu^\star=\inf\{\mu(P)|P\in\mathcal C\}$. 
Since the objective function~$\mu$ is continuous and the feasible
set~$\mathcal C$ is compact, there is at least one optimal transition
matrix $P^\star$ such that $\mu(P^\star)=\mu^\star$.
Let $\overline{P}$ denote the average over the orbit of $P^\star$
under $\mbox{Aut}(\mathcal G)$
\[
\overline{P}= \frac{1}{|\mbox{Aut}(\mathcal G)|} 
\sum_{\sigma\in\mbox{Aut}(\mathcal G)} \sigma(P^\star).
\]
This matrix is feasible because each $\sigma(P^\star)$ is feasible and
the feasible set is convex.
By construction, it is also invariant under the actions of
$\mbox{Aut}(\mathcal G)$. 
Moreover, using the convexity of~$\mu$, we have
$\mu(\overline{P})\leq\mu(P^\star)$. 
It follows that $\overline{P}\in\mathcal F$ and 
$\mu(\overline{P})=\mu^\star$.
\end{proof}

As a result of Theorem~\ref{t-fmmc-fixed-point}, we can replace the
constraint $P\in\mathcal C$ by $P\in\mathcal F$ in the FMMC problem
and get the same optimal value.
In the fixed-point subset $\mathcal F$, the transition
probabilities on the edges within an orbit must be the same. 
So we have the following corollaries:
\begin{corollary}\label{c-numprobs}
The number of distinct edge transition probabilities we need to
consider in the FMMC problem is at most equal to the number of 
orbits of $\,\mathcal{E}$ under $\mbox{Aut}(\mathcal G)$.
\end{corollary}
\begin{corollary}\label{c-edge-transitive}
If $\mathcal{G}$ is edge-transitive, then all the edge transition
probabilities can be assigned the same value.
\end{corollary}
Note that the holding probability at the vertices can always be
eliminated using $P_{ii} = 1- \sum_j P_{ij}$.
So it suffices to only consider the edge transition probabilities. 

\subsection{Formulation with reduced number of variables}
\label{s-reduce-variables}

From the results of the previous section, we can reduce the number of
optimization variables in the FMMC problem from the number of edges to
the number of edge orbits under the automorphism group. 
Here we give an explicit parametrization of the FMMC problem with
the reduced number of variables.
This parametrization is also the precise characterization of the
fixed-point subset $\mathcal F$. 

Recall that the \emph{adjacency matrix} of a graph with~$n$ vertices
is a $n\times n$ matrix~$A$ whose entries are given by $A_{ij}=1$ if
$\{i,j\}\in\mathcal{E}$ and $A_{ij}=0$ otherwise.
Let $\nu_i$ be the valency (degree) of vertex~$i$.
The \emph{Laplacian matrix} of the graph is given by 
$L=\Diag(\nu_1,\nu_2,\ldots,\nu_n)-A$, where $\Diag(\nu)$ denotes a
diagonal matrix with the vector $\nu$ as its diagonal. 
Extensive account of the Laplacian matrix and its use in algebraic
graph theory are provided in, \eg, \cite{Mer:94,FanChung,GoR:01}.

Suppose that there are $N$ orbits of edges under the action of 
$\mbox{Aut}(\mathcal G)$.
For each orbit, we define an orbit graph 
$\mathcal{G}_k=(\mathcal{V},\mathcal{E}_k)$, where $\mathcal{E}_k$ is
the set of edges in the $k$th orbit.
Note that the orbit graphs are disconnected (there are disconnected
vertices) if the original graph is not edge-transitive. 
Let $L_k$ be the Laplacian matrix of $\mathcal{G}_k$.
Note that the diagonal entries $(L_k)_{ii}$ equals the valency of
node~$i$ in $\mathcal{G}_k$ (which is zero if vertex~$i$ is
disconnected with all other vertices in $\mathcal{G}_k$).  

By Corollary~\ref{c-numprobs}, we can assign the same transition
probability on all the edges in the $k$-th orbit. 
Denote this transition probability by $p_k$ and let
$p=(p_1,\ldots,p_N)$.
Then the transition probability matrix can be written as 
\begin{equation}\label{e-prob-laplacian}
P(p) = I - \sum_{k=1}^N p_k L_k.
\end{equation}
This parametrization of the transition probability matrix  
automatically satisfies the constraints $P=P^T$, $P\ones=\ones$, 
and $P_{ij}=0$ for $\{i,j\}\in\mathcal E$. 
The entry-wise nonnegative constraint $P\geq 0$ now translates into
\begin{eqnarray*}
&& p_k \geq 0, \qquad k=1,\ldots,N \\
&& \sum_{k=1}^N (L_k)_{ii}\; p_k \leq 1, \qquad i=1,\ldots,n
\end{eqnarray*}
where the first set of constraints are for the off-diagonal entries
of~$P$, and the second set of constraints are for the diagonal entries
of~$P$. 

It can be verified that the parametrization~(\ref{e-prob-laplacian}),
together with the above inequality constraints, is the precise
characterization of the fixed-point subset~$\mathcal F$.
Therefore we can explicitly write the FMMC problem restricted to the
fixed-point subset as
\begin{equation}\label{e-fmmc-orbits}
\begin{array}{ll}
\mbox{minimize}   & \mu\left(I-\sum_{k=1}^N p_k L_k\right) \\[1ex]
\mbox{subject to} & p_k \geq 0, \quad k=1,\ldots,N\\[1ex]
& \sum_{k=1}^N (L_k)_{ii}\, p_k \leq 1, \quad i=1,\ldots,n.
\end{array}
\end{equation}
Later in this paper, we will also need the corresponding SDP
formulation 
\begin{equation}\label{e-fmmc-orbits-sdp}
\begin{array}{ll}
\mbox{minimize}   & s \\[0.5ex]
\mbox{subject to} & -s I \preceq I - \sum_{k=1}^N p_k L_k 
  - (1/n)\ones\ones^T \preceq s I\\[1ex]
& p_k \geq 0, \quad k=1,\ldots,N\\[1ex]
& \sum_{k=1}^N (L_k)_{ii}\, p_k \leq 1, \quad i=1,\ldots,n.
\end{array}
\end{equation}

\section{Some analytic results}
\label{s-analytic}

For some special classes of graphs, the FMMC problem can be
considerably simplified and often solved by only exploiting symmetry. 
In this section, we give some analytic results for the FMMC problem on
edge-transitive graphs, Cartesian product of simple graphs, and 
distance-transitive graphs (a subclass of
edge-transitive graphs). 
The optimal solution is often expressed in terms of the eigenvalues
of the adjacency matrix or the Laplacian matrix of the graph. 
It is interesting to notice that even for such highly structured class
of graphs, neither the maximum-degree nor the
Metropolis-Hastings heuristics discussed in~\cite{BDX:04} give the
optimal solution.
Throughout, we use $\alpha^\star$ to denote the common edge weight of the fastest mixing chain and $\mu^\star$ to denote the optimal SLEM. 

\subsection{FMMC on edge-transitive graphs}

\begin{theorem}\label{t-edge-transitive}
Suppose the graph $\mathcal G$ is edge-transitive, and let~$\alpha$ be
the transition probability assigned on all the edges.
Then the optimal solution of the FMMC problem is
\begin{eqnarray}
\alpha^\star 
&=& \min \left\{ \frac{1}{\nu_\mathrm{max}},\; 
\frac{2}{\lambda_1(L)+\lambda_{n-1}(L)} \right\}
\label{e-et-alpha}\\[1ex] 
\mu^\star 
&=& \max \left\{ 1-\frac{\lambda_{n-1}(L)}{\nu_\mathrm{max}},\; 
\frac{\lambda_1(L)-\lambda_{n-1}(L)}{\lambda_1(L)+\lambda_{n-1}(L)}
\right\}, 
\label{e-et-mu}
\end{eqnarray}
where $\nu_\mathrm{max}=\max_{i\in\mathcal{V}} \nu_i$ 
is the maximum valency of the vertices in the graph, and $L$ is the Laplacian matrix defined in~\S\ref{s-reduce-variables}.  
\end{theorem}

\begin{proof}
By definition of an edge-transitive graph, there is a single orbit of
edges under the actions of its automorphism group.  
Therefore we can assign the same transition probability $\alpha$ on
all the edges in the graph (Corollary~\ref{c-edge-transitive}), 
and the parametrization~(\ref{e-prob-laplacian})
becomes $P = I - \alpha L$. 
So we have
\[
\lambda_i(P) = 1-\alpha \lambda_{n+1-i}(L), \qquad i=1,\ldots,n 
\]
and the SLEM
\begin{eqnarray*}
\mu(P) &=& \max\{\lambda_2(P),\;-\lambda_{n}(P)\}\\
&=&\max\{ 1-\alpha \lambda_{n-1}(L), \; \alpha \lambda_1(L)-1 \}.
\end{eqnarray*}
To minimize $\mu(P)$, we let 
$1-\alpha \lambda_{n-1}(L)=\alpha \lambda_1(L)-1 $ and 
get $\alpha=2/(\lambda_{n-1}(L)+\lambda_{n-1}(L))$.  
But the nonnegativity constraint $P\geq 0$ requires that 
the transition probability must also satisfy 
$0<\alpha\leq 1/\nu_\mathrm{max}$. 
Combining these two conditions gives the optimal 
solution~(\ref{e-et-alpha}) and~(\ref{e-et-mu}). 
\end{proof}

We give two examples of FMMC on edge-transitive graphs. 

\subsubsection{Cycles}

\begin{figure}
\centering
\includegraphics[width=0.2\textwidth]{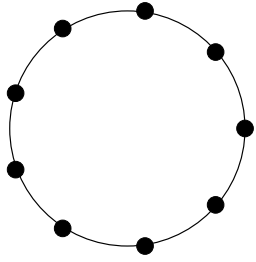}
\caption{The cycle graph $C_n$ with $n=9$.}
\label{f-cycle}
\end{figure}

The first example is the cycle graph $C_n$; see Figure~\ref{f-cycle}. 
The Laplacian matrix is 
\[
L=\left[\begin{array}{rrrrrr}
 2 & -1 &  0 & \cdots & 0 & -1 \\
-1 &  2 & -1 & \cdots & 0 &  0 \\
 0 & -1 &  2 & \cdots & 0 &  0 \\
\vdots & \vdots & \vdots & \ddots &\vdots & \vdots \\  
 0 &  0 &  0 & \cdots & 2 & -1 \\
-1 &  0 &  0 & \cdots & -1 & 2 \end{array} \right]
\]
which has eigenvalues
\[
 2 - 2\cos\frac{2k\pi}{n}, \qquad k=1,\ldots,n.
\]
The two extreme eigenvalues are
\[
\lambda_1(L)= 2 -  2 \cos \frac{2\lfloor n/2\rfloor\pi}{n}, \qquad
\lambda_{n-1}(L) = 2 - 2 \cos \frac{2\pi}{n}
\]
where $\lfloor n/2\rfloor$ denotes the largest integer that is
no larger than $n/2$, which is $n/2$ for $n$ even 
or $(n-1)/2$ for $n$ odd. 
By Theorem~\ref{t-edge-transitive}, 
the optimal solution to the FMMC problem is 
\begin{eqnarray}
\alpha^{\star}
&=&\frac{1}{2-\cos\frac{2\pi}{n}-\cos\frac{2\lfloor 
   n/2\rfloor\pi}{n}} \label{e-cycle-prob}\\[2ex]
\mu^{\star}
&=&\frac{\cos\frac{2\pi}{n}-\cos\frac{2\lfloor n/2\rfloor\pi}{n}}
   {2-\cos\frac{2\pi}{n}-\cos\frac{2\lfloor n/2\rfloor\pi}{n}}.
   \label{e-cycle-slem}
\end{eqnarray}
When $n\to\infty$, the transition probability $\alpha^{\star}\to 1/2$
and the SLEM $\mu^{\star}\to 1-2\pi^2/n^2$.

\subsubsection{Complete bipartite graphs}
\begin{figure}
\centering
\psfrag{u}[cr]{$u\;$}
\psfrag{v}[cl]{$\;v$}
\psfrag{x}[cr]{$x\;$}
\psfrag{y}[cl]{$\;y$}
\includegraphics[width=0.22\textwidth]{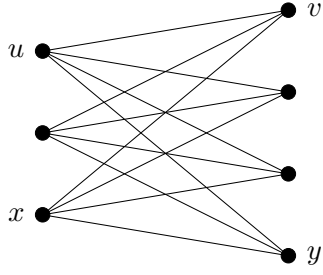}
\caption{The complete bipartite graph $K_{m,n}$ with $m=3$ and $n=4$.}
\label{f-bipartite}
\end{figure}

The complete bipartite graph, denoted $K_{m,n}$, has two subsets of
vertices with cardinalities~$m$ and~$n$ respectively.
Each vertex in a subset is connected to all the vertices in the other
subset,  and is not connected to any of the vertices in its own
subset; see Figure~\ref{f-bipartite}.  
Without loss of generality, assume $m\leq n$.
So the maximum degree is~$\nu_\mathrm{max}=n$. 
This graph is edge-transitive but not vertex-transitive.  
The Laplacian matrix of this graph is
\[
L=\left[\begin{array}{cc} n I_m & -\ones_{m\times n} \\
 -\ones_{n\times m} & m I_n \end{array}\right]
\]
where $I_m$ denotes the~$m$ by~$m$ identity matrix, and 
$\ones_{m\times n}$ denotes the~$m$ by~$n$ matrix whose components are
all ones.  
For $n\geq m\geq 2$, this matrix has four distinct eigenvalues
$m+n$, $n$, $m$ and $0$, with multiplicities $1$, $m-1$, $n-1$ and
$1$, respectively (see~\S\ref{s-bipartite-diag}). 
By Theorem~\ref{t-edge-transitive}, the optimal transition probability
on the edges and the corresponding SLEM are
\begin{eqnarray}
\alpha^{\star} &=& \min \left\{ \frac{1}{n}, \;\frac{2}{n+2m} \right\} 
   \label{e-bipartite-alpha} \\[1ex]
\mu^{\star} &=& \max \left\{ \frac{n-m}{n},\;\frac{n}{n+2m} \right\}.
   \label{e-bipartite-mu}
\end{eqnarray}

\subsection{Cartesian product of graphs}

Many graphs we consider can be constructed
by taking \emph{Cartesian product} of simpler graphs.
The Cartesian product of two graphs
$\mathcal{G}_1=(\mathcal{V}_1,\mathcal{E}_1)$ and 
$\mathcal{G}_2=(\mathcal{V}_2,\mathcal{E}_2)$ is a graph with vertex
set $\mathcal{V}_1 \times \mathcal{V}_2$, where two vertices 
$(u_1,u_2)$ and $(v_1, v_2)$ are connected by
an edge if and only if $u_1=v_1$ and $\{u_2,v_2\}\in \mathcal{E}_2$,
or $u_2=v_2$ and $\{u_1,v_1\}\in \mathcal{E}_1$. 
Let $\mathcal G_1\oplus\mathcal G_2$ denote this Cartesian product.
Its Laplacian matrix is given by
\begin{equation}\label{e-Cartesian-laplacian}
L_{\mathcal G_1\oplus\mathcal G_2} = L_{\mathcal G_1}\otimes 
I_{|\mathcal V_1|} + I_{|\mathcal V_2|}\otimes L_{\mathcal G_2}
\end{equation}
where 
$\otimes$ denotes the matrix Kronecker product (\cite{Gra:81}). 
The eigenvalues of $L_{\mathcal G_1\oplus\mathcal G_2}$ are given by 
\begin{equation}\label{e-Cartesian-eigs}
\lambda_i(L_{\mathcal G_1}) + \lambda_j(L_{\mathcal G_2}), 
\qquad i=1,\ldots,|\mathcal{V}_1|, \qquad j=1,\ldots,|\mathcal{V}_2|
\end{equation}
where each eigenvalue is obtained as many times as its multiplicity
(\eg, \cite{Moh:97}). 
The adjacency matrix of the Cartesian product of graphs also has
similar properties, which we will use later for distance-transitive
graphs. 
Detailed background on spectral graph theory can be found in, \eg,
\cite{Biggs,DCS:80,FanChung,GoR:01}. 

Combining Theorem~\ref{t-edge-transitive} and the above expression for
eigenvalues, we can easily obtain solutions to the FMMC problem on  
graphs formed by taking Cartesian product of simple graphs. 

\subsubsection{Two-dimensional meshes}
\label{s-twodimgrid}

\begin{figure}
\centering
\includegraphics[width=0.22\textwidth]{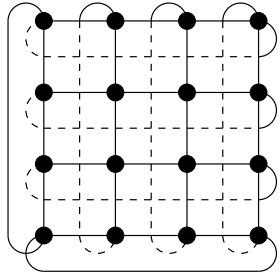}
\caption{The two-dimensional mesh with wraparounds $M_n$ with $n=4$.}
\label{f-mesh}
\end{figure}

Here we consider the two-dimensional mesh with wraparounds at two ends
of each row and column, see Figure~\ref{f-mesh}. 
It is the Cartesian product of two copies of $C_n$. 
We write it as $M_n=C_n\oplus C_n$. 
By equation~(\ref{e-Cartesian-eigs}), its Laplacian matrix has eigenvalues  
\[
4-2\cos\frac{2i\pi}{n}-2\cos\frac{2j\pi}{n},
\qquad i,~j=1,\ldots,n.
\]
By Theorem~\ref{t-edge-transitive}, we obtain
the optimal transition probability
\[
\alpha^{\star}=\frac{1}{3-2\cos\frac{2\lfloor n/2\rfloor\pi}{n}
-\cos\frac{2\pi}{n}}
\]
and the smallest SLEM
\[
\mu^{\star}=\frac{1-2\cos\frac{2\lfloor n/2\rfloor\pi}{n}+
 \cos\frac{2\pi}{n}}{3-2\cos\frac{2\lfloor n/2\rfloor\pi}{n}
-\cos\frac{2\pi}{n}}
\]
When $n\to\infty$, the transition probability $\alpha^{\star}\to 1/4$
and the SLEM $\mu^{\star}\to 1-\pi^2/n^2$.

\subsubsection{Hypercubes}
\label{s-hypercubes}
\begin{figure}
\centering
\includegraphics[width=0.45\textwidth]{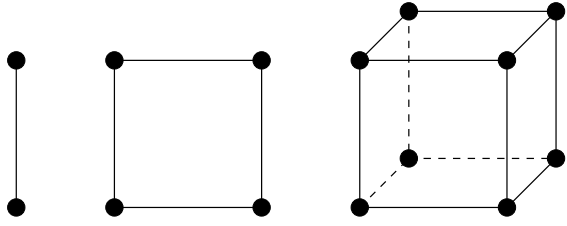}
\caption{The hypercubes $Q_1,~Q_2$ and
$Q_3$.} 
\label{f-hypercube}
\end{figure}

The $d$-dimensional hypercube, denoted $Q_d$, has $2^d$ vertices, each
labeled with a binary word with length~$d$. 
Two vertices are connected by an edge if their words differ in exactly
one component (see Figure~\ref{f-hypercube}). 
This graph is isomorphic to the Cartesian product of $d$ copies of
$K_2$, the complete graph with two vertices. 
The Laplacian of $K_2$ is
\[
L_{K_2} = \left[ \begin{array}{rr} 1 & -1\\ -1 & 1\end{array}\right],
\]
whose two eigenvalues are~$0$ and~$2$.
The one-dimensional hypercube $Q_1$ is just $K_2$.
Higher dimensional hypercubes are defined recursively:
\[
Q_{k+1} = Q_k \oplus K_2, \qquad k=1,2,\ldots.
\] 
By equation~(\ref{e-Cartesian-laplacian}), their Laplacian matrices are
\[
L_{Q_{k+1}} = L_{Q_k} \otimes I_2 + I_{2^k} \otimes L_{K_2}, \qquad
k=1,2,\ldots .
\]
Using equation~(\ref{e-Cartesian-eigs}) recursively, 
the Laplacian of $Q_d$ has eigenvalues $2k$, $k=0,1,\ldots,d$, each
with multiplicity $\left({d\atop k}\right)$.
The FMMC is achieved for:
\[
\alpha^{\star}=\frac{1}{d+1}, \qquad
\mu^{\star}=\frac{d-1}{d+1}.
\]
This solution has also been worked out, for example, in~\cite{Moh:97}. 

\subsection{FMMC on distance-transitive graphs}

Distance-transitive graphs have been studied extensively in the
literature (see, \eg, \cite{BCN:89}).  In particular, they are both
edge- and vertex-transitive.  In previous examples, the cycles and the
hypercubes are actually distance-transitive graphs; so are the
bipartite graphs when the two parties have equal number of vertices.

In a distance-transitive graph, all vertices have the same valency,
which we denote by~$\nu$.
The Laplacian matrix can be written as $L = \nu I - A$, with $A$ being
the adjacency matrix. 
Therefore 
\[
\lambda_i(L) = \nu - \lambda_{n+1-i}(A), \qquad i=1,\ldots,n.
\]
We can substitute the above equation in  
equations~(\ref{e-et-alpha}) and~(\ref{e-et-mu}) to obtain the optimal
solution in terms of $\lambda_2(A)$ and $\lambda_n(A)$.

Since distance-transitive graphs usually have very large automorphism 
groups, the eigenvalues of the adjacency matrix~$A$ (and the
Laplacian~$L$) often have very high multiplicities. 
But to solve the FMMC problem, we only need to know the distinct
eigenvalues; actually, only $\lambda_2(A)$ and $\lambda_n(A)$ would
suffice.  
In this regard, it is more convenient to use a much smaller matrix, 
the \emph{intersection matrix}, which has all the distinct eigenvalues
of the adjacency matrix. 

Let $D$ be the \emph{diameter} of the graph. 
For a nonnegative integer~$k\leq D$, choose any two vertices~$u$
and~$v$ such that their distance satisfies $\delta(u,v)=k$.
Let $a_k$, $b_k$ and $c_k$ be the number of vertices that 
are adjacent to~$u$ and whose distance from~$v$ are $k$, $k+1$ and
$k-1$, respectively. That is,
\begin{eqnarray*}
a_k &=& |\{ w\in\mathcal V ~|~ \delta(u,w)=1, \; \delta(w,v)=k\}|\\ 
b_k &=& |\{ w\in\mathcal V ~|~ \delta(u,w)=1, \; \delta(w,v)=k+1\}|\\ 
c_k &=& |\{ w\in\mathcal V ~|~ \delta(u,w)=1, \; \delta(w,v)=k-1\}|. 
\end{eqnarray*}
For distance-transitive graphs, these numbers are independent
of the particular pair of vertices~$u$ and~$v$ chosen. 
Clearly, we have $a_0=0$, $b_0=\nu$ and $c_1=1$. 
The \emph{intersection matrix} $B$ is the following tridiagonal 
$(D+1)\times(D+1)$ matrix
\[
B = \left[ \begin{array}{ccccc}
a_0 & b_0 &        &        &        \\
c_1 & a_1 & b_1    &        &        \\
    & c_2 & a_2 & \ddots &        \\
    &     & \ddots & \ddots & b_{D-1}\\
    &     &        & c_D    & a_D 
\end{array} \right]. 
\]

Denote the eigenvalues of the intersection matrix, in decreasing
order, as $\eta_0, ~\eta_1, ~\ldots, ~\eta_D$.  These are precisely
the $(D+1)$ distinct eigenvalues of the adjacency matrix~$A$ (see,
\eg, \cite{Biggs}).  In particular, we have
\[
\lambda_1(A) = \eta_0 = \nu, \qquad \lambda_2(A) = \eta_1, \qquad 
\lambda_n(A) = \eta_D. 
\]
The following corollary is a direct consequence of
Theorem~\ref{t-edge-transitive}.
\begin{corollary}\label{c-distance-transitive}
The optimal solution of the FMMC problem on a distance-transitive
graph is
\begin{eqnarray}
\alpha^\star &=& \min\left\{ \frac{1}{\nu},\; 
  \frac{2}{2\nu -(\eta_1 + \eta_D)} \right\} \label{e-dist-tran-alpha}\\
\mu^\star &=& \max\left\{ \frac{\eta_1}{\nu}, \;
  \frac{\eta_1 - \eta_D}{2\nu - (\eta_1 + \eta_D)} \right\} .
\label{e-dist-tran-mu}
\end{eqnarray}
\end{corollary}

Next we give solutions for the FMMC problem on several families of
distance-transitive graphs.

\subsubsection{Complete graphs}
The case of the complete graph with $n$ vertices, usually called
$K_n$, is very simple.
It is distance-transitive, with diameter $D=1$ and valency $\nu=n-1$.
The intersection matrix is
\[
B=\left[
\begin{array}{cc}
 0 & n-1 \\ 1 & n-2
\end{array}
\right],
\]
with eigenvalues $\eta_0=n-1$, $\eta_1=-1$.  
Using equations~(\ref{e-dist-tran-alpha}) and~(\ref{e-dist-tran-mu}), 
the optimal parameters are
\[
\alpha^{\star}=\frac{1}{n}, \qquad \mu^{\star}=0.
\]
The associated matrix $P=(1/n)\ones\ones^T$ has one
eigenvalue equal to 1, and the remaining $n-1$ eigenvalues vanish.
Such Markov chains achieve perfect mixing after just one step,
regardless of the value of $n$.

\subsubsection{Petersen graph}

\begin{figure}
\begin{center}
\includegraphics[width=3.5cm]{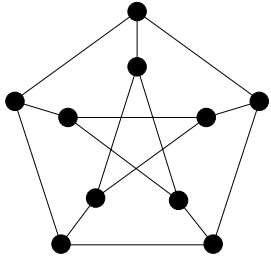}
\caption{The Petersen graph.}
\label{fig:petersen}
\end{center}
\end{figure}

The Petersen graph, shown in Figure~\ref{fig:petersen}, is a well-known
distance-transitive graph with 10 vertices and 15 edges. The diameter
of the graph is $D=2$, and the intersection matrix is 
\[
B=\left[
\begin{array}{ccc}
0 & 3 & 0 \\
1 & 0 & 2 \\
0 & 1 & 2 
\end{array}
\right]
\]
with eigenvalues $\eta_0 =3$, $\eta_1 = 1$ and $\eta_2 = -2$. 
Applying the formula~(\ref{e-dist-tran-alpha})
and~(\ref{e-dist-tran-mu}), we obtain
\[
\alpha^\star = \frac{2}{7}, \qquad \mu^\star = \frac{3}{7}.
\]

\subsubsection{Hamming graphs}

The Hamming graphs, denoted $H(d,n)$, have vertices labeled by
elements in the Cartesian product $\{1,\ldots,n\}^d$, with two
vertices being adjacent if they differ in exactly one component.
By the definition, it is clear that Hamming graphs are isomorphic to
the Cartesian product of $d$ copies of the complete graph $K_n$.
Hamming graphs are distance-transitive, with diameter $D=d$ and
valency $\nu=d\,(n-1)$.  Their eigenvalues are given by 
$\eta_k = d \, (n-1) - k n$ for $k = 0, \ldots, d$.
These can be obtained using an equation for eigenvalues of adjacency
matrices, similar to~(\ref{e-Cartesian-eigs}), 
with the eigenvalues of~$K_n$ being $n-1$ and~$-1$.  
Therefore the FMMC has parameters:
\begin{eqnarray*}
\alpha^\star 
&=& \min\left\{ \frac{1}{d \,(n-1)}, \; \frac{2}{n \, (d+1)} 
\right\}\\[1ex]
\mu^\star 
&=& \max\left\{ 1-\frac{n}{d(n-1)}, \; \frac{d-1}{d+1} \right\}. 
\end{eqnarray*}
We note that hypercubes (see \S\ref{s-hypercubes}) are special Hamming
graphs with $n=2$. 

\subsubsection{Johnson graphs}
The Johnson graph $J(n,q)$ (for $1 \leq q \leq n/2$) is defined as
follows: the vertices are the $q$-element subsets of $\{1,\ldots,n\}$,
with two vertices being connected with an edge if and only if the
subsets differ exactly by one element. It is a distance-transitive
graph, with $n \choose q$ vertices and $\frac{1}{2} q \,(n-q) {n
\choose q}$ edges.  It has valency $\nu=q\,(n-q)$ and diameter $D=q$.
The eigenvalues of the intersection matrix can be computed
analytically and they are:
\[
\eta_k = q\,(n-q) + k\,(k-n-1), \qquad k = 0, \ldots, q. 
\]
Therefore, by Corollary~\ref{c-distance-transitive}, we obtain
the optimal transition probability
\[
\alpha^\star = 
\min\left\{\frac{1}{q \,(n-q)}, ~\frac{2}{q n+n+q-q^2}\right\}
\]
and the smallest SLEM
\[
\mu^\star=\max\left\{1-\frac{n}{q(n-q)},~1-\frac{2n}{q n+n+q-q^2}\right\}.
\]

\section{FMMC on orbit graphs}
\label{s-fmmc-orbit}
For graphs with large automorphism groups, the eigenvalues of the
transition probability matrix often have very high multiplicities.  To
solve the FMMC problem, it suffices to work with only the distinct
eigenvalues without consideration of their multiplicities.  This is
exactly what the intersection matrix does for distance-transitive
graphs.  In this section we develop similar tools for more general
graphs.  More specifically, we show how to construct an orbit chain
which is much smaller in size than the original Markov chain, but
contains all its distinct eigenvalues (with much fewer
multiplicities).  The FMMC on the original graph can be found by
solving a much smaller problem on the orbit chain.

\subsection{Orbit theory}
\label{s-orbit-theory}
Here we review the orbit theory developed in \cite{BDPX:05}. 
Let $P$ be a symmetric Markov chain on the graph 
$\mathcal G=(\mathcal V, \mathcal E)$, and $H$ be a group of automorphisms of the graph. 
Often, it is a subgroup of the full automorphism group
$\mbox{Aut}(\mathcal G)$.
The vertex set $\mathcal V$ partitions into orbits
$O_v = \{hv:h\in H\}$.
For notational convenience, in this section we use $P(v,u)$, 
for $v,u\in\mathcal V$,  
to denote entries of the transition probability matrix. 
We define the \emph{orbit chain} by
specifying the transition probabilities between orbits
\begin{equation}\label{e-orbit-transition}
P_H(O_v, O_u) = P(v,O_u) = \sum_{u'\in O_u} P(v,u').
\end{equation}
This transition probability is independent of which 
$v\in O(v)$ is chosen, so it is well defined 
and the lumped orbit chain is indeed Markov.

The orbit chain is in general no longer symmetric, 
but it is always reversible. 
Let $\pi(i)$, $i\in\mathcal V$, be the stationary distribution of the
original Markov chain.
Then the stationary distribution on the orbit chain is obtained as
\begin{equation}\label{e-orbit-distribution}
\pi_H(O_v)=\sum_{i\in O_v} \pi(i) .
\end{equation}
It can be verified that 
\begin{equation}\label{e-detailed-balance}
\pi_H(O_v) P_H(O_v,O_u) = \pi_H(O_u) P_H(O_u,O_v),
\end{equation}
which is the detailed balance condition to test reversibility. 

The following is a summary of the orbit theory we developed in
\cite{BDPX:05}, which relate the eigenvalues and eigenvectors of
the orbit chain $P_H$ to the eigenvalues and eigenvectors of the
original chain $P$.
\begin{itemize}
\item \emph{Lifting} (\cite[\S3.1]{BDPX:05}).
If $\bar \lambda$ is an eigenvalue of $P_H$ with associated
eigenvector $\bar f$, then $\bar \lambda$ is an eigenvalue of $P$ with
$H$-invariant eigenfunction $f(v)=\bar f(O_v)$.
Conversely, every $H$-invariant eigenfunction appears uniquely from
this construction.
\item \emph{Projection} (\cite[\S3.2]{BDPX:05}).
Let $\lambda$ be an eigenvalue of $P$ with eigenvector $f$. 
Define $\bar f(O_v)=\sum_{h\in H} f(h^{-1}(v))$.
Then $\lambda$ appears as an eigenvalue of $P_H$, with eigenvector
$\bar f$, if either of the following
two conditions holds:
\begin{itemize}
\item[(a)] $H$ has a fixed point $v^*$ and $f(v^*)\neq 0$.
\item[(b)] $f$ is nonzero at a vertex $v^*$ in an 
$\mbox{Aut}(\mathcal G)$-orbit which contains a fixed point of~$H$.
\end{itemize}
\end{itemize}

Equipped with this orbit theory, we would like to construct one or
multiple orbit chains that retain all the eigenvalues of the original
chain.  Ideally the orbit chains are much smaller in size than the
original chain, with the eigenvalues having much fewer multiplicities.
The following theorem (Theorem 3.7 in \cite{BDPX:05}) gives sufficient
conditions that guarantee that the orbit chain(s) attain all the
eigenvalues of the original chain.

\begin{theorem}\label{t-full-eig-proj}
Suppose that $\mathcal V=O_1\cup\ldots\cup O_K$ is a disjoint union of
the orbits under $\mbox{Aut}(\mathcal G)$.
Let $H_i$ be the subgroup of $\mbox{Aut}(\mathcal G)$ that has a fixed
point in $O_i$. 
Then all eigenvalues of $P$ occur among the eigenvalues of 
$\{P_{H_i}\}_{i=1}^K$.
Further, every eigenvector of $P$ occurs by lifting an eigenvector of
some $P_{H_i}$.
\end{theorem}

Observe that if $H\subseteq G\subseteq \mbox{Aut}(\mathcal G)$, then
the eigenvalues of $P_H$ contain all eigenvalues of $P_G$.  
This allows disregarding some of the $H_i$ in
Theorem~\ref{t-full-eig-proj}. 
In particular, it is possible to construct a single orbit chain that
contains all eigenvalues of the original chain. 
Therefore we have

\begin{corollary} \label{c-full-eig-proj}
Suppose that $\mathcal V=O_1\cup\ldots\cup O_k$ is a disjoint union of
the orbits under $\mbox{Aut}(\mathcal G)$, and
$H$ is a subgroup of $\mbox{Aut}(\mathcal G)$.
If $H$ has a fixed point in every $O_i$, 
then all distinct eigenvalues of $P$ occur among the eigenvalues of 
$P_H$.
\end{corollary}

\paragraph{Remarks.}
To find $H$ in the above corollary, we can just compute the corresponding
stabilizer, \ie, compute the largest subgroup of $\mbox{Aut}(\mathcal G)$
that fixes one point in each orbit. 
Note that the $H$ promised by the corollary may be trivial in some cases; 
see the example in~\S\ref{s-graham-example}.
\vspace{2ex}

\begin{figure}[t]
\centering
\subfigure[\normalsize Orbit chain under $S_m\times S_n$.]{
\label{f-bipartite-orbit:a}
\begin{minipage}[t]{0.49\textwidth}
\centering
\psfrag{mp}[tc]{$mp$}
\psfrag{np}[bc]{$np$}
\psfrag{u}[cr]{$O_u$}
\psfrag{v}[cl]{$O_v$}
\includegraphics[width=0.5\textwidth]{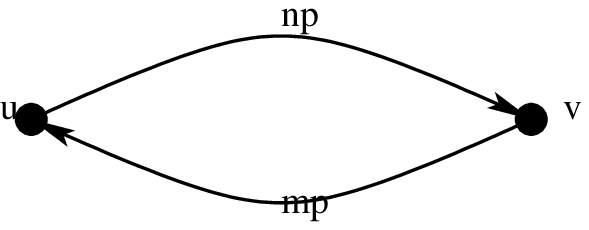}
\end{minipage}}%
\subfigure[\normalsize Orbit chain under $S_{m-1}\times S_n$.]{
\label{f-bipartite-orbit:b}
\begin{minipage}[t]{0.49\textwidth}
\centering
\psfrag{np}[bc]{$np$}
\psfrag{pbr}[br]{$p$}
\psfrag{n-1p}[tl]{$np$}
\psfrag{m-1pc}[bc]{$(m\!-\!1)p$}
\psfrag{x}[cr]{$x$}
\psfrag{u}[cr]{$O_u$}
\psfrag{v}[cl]{$O_v$}
\includegraphics[width=0.5\textwidth]{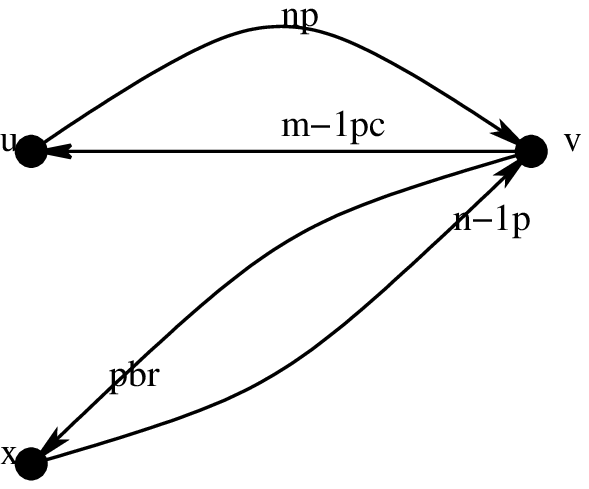}
\end{minipage}}
\newline\\[2ex]
\subfigure[\normalsize Orbit chain under $S_m\times S_{n-1}$.]{
\label{f-bipartite-orbit:c}
\begin{minipage}[t]{0.49\textwidth}
\centering
\psfrag{mp}[bc]{$mp$}
\psfrag{pbl}[bl]{$\!~p$}
\psfrag{m-1p}[tr]{$mp$}
\psfrag{n-1pc}[bc]{$(n\!-\!1)p$}
\psfrag{y}[cl]{$y$}
\psfrag{u}[cr]{$O_u$}
\psfrag{v}[cl]{$O_v$}
\includegraphics[width=0.5\textwidth]{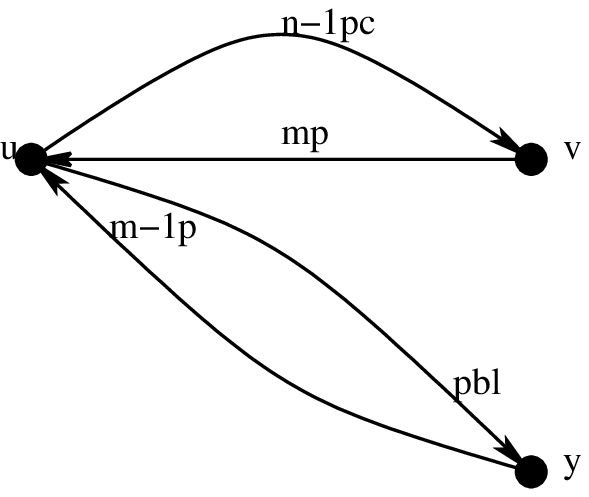}
\end{minipage}}%
\subfigure[\normalsize Orbit chain under $S_{m-1}\times S_{n-1}$.]{
\label{f-bipartite-orbit:d}
\begin{minipage}[t]{0.49\textwidth}
\centering
\psfrag{p}[bc]{$p$}
\psfrag{pbl}[bl]{$\!~p$}
\psfrag{pbr}[br]{$p$}
\psfrag{n-1p}[tl]{$\!(n\!-\!1)p$}
\psfrag{m-1p}[tr]{$(m\!-\!1)p\!$}
\psfrag{n-1pc}[bc]{$(n\!-\!1)p$}
\psfrag{m-1pc}[tc]{$(m\!-\!1)p$}
\psfrag{x}[cr]{$x$}
\psfrag{y}[cl]{$y$}
\psfrag{u}[cr]{$O_u$}
\psfrag{v}[cl]{$O_v$}
\includegraphics[width=0.5\textwidth]{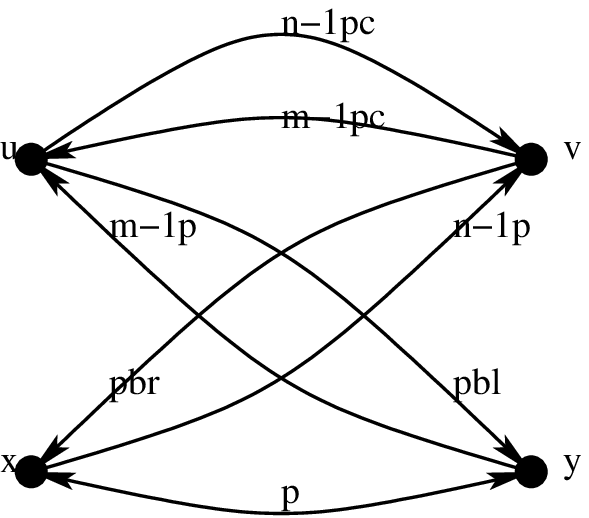}
\end{minipage}} 
\caption{Orbit chains of $K_{m,n}$ under different automorphism groups.
 The vertices labeled $O_u$ and $O_v$ are orbits of vertices $u$ and
  $v$ (labeled in Figure~\ref{f-bipartite}) under corresponding
  actions. The vertices labeled~$x$ and~$y$ are fixed points.} 
\label{f-bipartite-orbit}
\end{figure}

We illustrate the orbit theory with the bipartite graph $K_{m,n}$ 
shown in Figure~\ref{f-bipartite}.  
It is easy to see that $\mbox{Aut}(K_{m,n})$ is the direct product of
two symmetric groups, namely $S_m \times S_n$, with each symmetric
group permuting one of the two subsets of vertices. 
This graph is edge-transitive.
So we assign the same transition probability~$p$ on all the edges. 

The orbit chains under four different subgroups of
$\mbox{Aut}(K_{m,n})$ are shown in Figure~\ref{f-bipartite-orbit}.
The transition probabilities between orbits are calculated using
equation~(\ref{e-orbit-transition}).
Since the transition probabilities are not symmetric, we represent
the orbit chains by directed graphs, with different transition 
probabilities labeled on opposite directions between two adjacent
vertices. 
The full automorphism group $\mbox{Aut}(K_{m,n})$ has two orbits of
vertices; see Figure~\ref{f-bipartite-orbit:a}. 
The orbit graphs under the subgroups $S_{m-1}\times S_n$
(Figure~\ref{f-bipartite-orbit:b}) and $S_m\times S_{n-1}$ 
(Figure~\ref{f-bipartite-orbit:c}) each contains a fixed point of the
two orbits under $\mbox{Aut}(K_{m,n})$.
By Theorem~\ref{t-full-eig-proj}, these two orbit chains contain
all the distinct eigenvalues of the original chain on $K_{m,n}$.
Alternatively, we can construct the orbit chain under the subgroup 
$S_{m-1}\times S_{n-1}$, shown in Figure~\ref{f-bipartite-orbit:d}.
This orbit chain contain a fixed point in both orbits under
$\mbox{Aut}(K_{m,n})$. 
By Corollary~\ref{t-full-eig-proj}, all distinct eigenvalues of
$K_{m,n}$ appear in this orbit chain. 
In particular, this shows that there are at most four distinct
eigenvalues in the original chain.

If we order the vertices in Figure~\ref{f-bipartite-orbit:d} as
$(x,y,O_u, O_v)$, then the transition probability matrix for this
orbit chain is
\[
P_H = \left[ \begin{array}{cccc}
1-np & p & 0 & (n-1)p \\
p & 1-mp & (m-1)p & 0\\
0 & p & 1-np & (n-1)p \\
p & 0 & (m-1)p & 1-mp \end{array} \right]
\]
where $H=S_{m-1}\times S_{n-1}$.
By equation~(\ref{e-orbit-distribution}), 
its stationary distribution is
\[
\pi_H = \left( \frac{1}{m+n},~\frac{1}{m+n},~\frac{m-1}{m+n},
        ~\frac{n-1}{m+n} \right).
\]

\subsection{Fastest mixing reversible Markov chain on orbit graph}
\label{s-fmmc-reversible}

Since in general the orbit chain is no longer symmetric, we cannot
directly use the convex optimization formulation~(\ref{e-fmmc})
or~(\ref{e-fmmc-sdp}) to minimize $\mu(P_H)$.  Fortunately, the
detailed balance condition~(\ref{e-detailed-balance}) leads to a
simple transformation that allow us to formulate the problem of
finding the fastest reversible Markov chain as a convex program
\cite{BDX:04}.
 
Suppose the orbit chain $P_H$ contains all distinct eigenvalues of the 
original chain. 
Let $\pi_H$ be the stationary distribution of the orbits, and 
let $\Pi=\Diag(\pi_H)$.
The detailed balance condition~(\ref{e-detailed-balance}) can be 
written as $\Pi P_H = P_H^T \Pi$,
which implies that the matrix $\Pi^{1/2}P_H\Pi^{-1/2}$ is symmetric 
(and of course, has the same eigenvalues as $P_H$). 
The eigenvector of $\Pi^{1/2}P_H\Pi^{-1/2}$ associated with the
maximum eigenvalue~$1$ is 
$q=(\sqrt{\pi_H(O_1)},\ldots,\sqrt{\pi_H(O_k)})$.  
The SLEM $\mu(P_H)$ equals the spectral norm of
$\Pi^{1/2}P_H\Pi^{-1/2}$ restricted to the orthogonal complement of
the subspace spanned by $q$. 
This can be written as
\[
\mu(P_H) = \| (I-qq^T) \Pi^{1/2} P_H \Pi^{-1/2} (I-qq^T) \|_2
      = \| \Pi^{1/2} P_H \Pi^{-1/2} - qq^T \|_2.
\]
Introducing a scalar variable $s$ to bound the above spectral norm,  
we can formulate the fastest mixing reversible Markov chain problem as 
an SDP
\begin{equation}\label{e-fmmc-rev}
\begin{array}{ll}
\mbox{minimize}  & s \\[0.5ex]
\mbox{subject to} 
& -sI \preceq \Pi^{1/2} P_H \Pi^{-1/2} - qq^T \preceq sI\\[0.5ex]
& P_H\geq 0,\quad P_H\ones=\ones,\quad\Pi P_H=P_H^T\Pi\\[0.5ex]  
& P_H(O,O') = 0, \quad  (O,O')\notin \mathcal{E}_H.
\end{array}
\end{equation}
The optimization variables are the matrix $P_H$ and scalar~$s$, and
problem data is given by the orbit graph and the stationary
distribution $\pi_H$.  Note that the reversibility constraint $\Pi
P_H=P_H^T\Pi$ can be dropped since it is always satisfied by the
construction of the orbit chain; see
equation~(\ref{e-detailed-balance}). By pre- and post-multiplying the
matrix inequality by $\Pi^{1/2}$, we can write then another
equivalent formulation:
\begin{equation}\label{e-fmmc-rev2}
\begin{array}{ll}
\mbox{minimize}  & s \\[0.5ex]
\mbox{subject to} 
& -s \Pi \preceq \Pi P_H - \pi_H \pi_H^T \preceq s \Pi\\[0.5ex]
& P_H\geq 0,\quad P_H\ones=\ones, \\[0.5ex]  
& P_H(O,O') = 0, \quad  (O,O')\notin \mathcal{E}_H.
\end{array}
\end{equation}
To solve the fastest mixing reversible Markov chain problem on the 
orbit graph, we need the following three steps.
\begin{enumerate}
\item
Conduct symmetry analysis on the original graph: identify the
automorphism graph $\mbox{Aut}(\mathcal G)$ and determine the number
of orbits of edges $N$. 
By Corollary~\ref{c-numprobs},
this is the number of transition probabilities we need to consider.
\item 
Find a group of automorphisms $H$ that satisfies the conditions in
Corollary~\ref{c-full-eig-proj}. 
Construct its orbit chain by computing the transition probabilities
using equation~(\ref{e-orbit-transition}), and compute the stationary
distribution using equation~(\ref{e-orbit-distribution}). 
Note that the entries of $P_H$ are multiples of the transition
probabilities on the original graph.
\item
Solve the fastest mixing reversible Markov chain
problem~(\ref{e-fmmc-rev}).
The optimal SLEM $\mu(P_H^\star)$ is also the optimal SLEM for the
original chain, and the optimal transition probabilities on the
original chain can be obtained by simple scaling of the optimal orbit
transition probabilities. 
\end{enumerate}

We have assumed a single orbit chain that contains all distinct
eigenvalues of the original chain.
Sometimes it is more convenient to use multiple orbit chains.
Let $P_{H_i}$, $i=1,\ldots,K$, be the collection of orbit chains
in Theorem~\ref{t-full-eig-proj}.
In this case we need to minimize $\max_i \mu(P_{H_i})$. 
This can be done by simply adding the set of constraints 
in~(\ref{e-fmmc-rev}) for every matrix $P_{H_i}$.
For example, for the complete bipartite graph $K_{m,n}$, instead of
using the single orbit chain in Figure~\ref{f-bipartite-orbit:d},  
we can use the two orbit chains in Figure~\ref{f-bipartite-orbit:b} 
and Figure~\ref{f-bipartite-orbit:c} together, with two sets of
constraints in the SDP~(\ref{e-fmmc-rev}).

\subsection{Examples}
We demonstrate the above computational procedure on orbit graphs 
with two examples: the graph $K_n$-$K_n$ and the complete binary
tree. 
Both examples will be revisited in~\S\ref{s-block-diag} using the
method of block diagonalization.

\subsubsection{The graph $K_n$-$K_n$}
\label{s-Kn-Kn}

\begin{figure}
\centering
\subfigure[\normalsize The graph $K_n$-$K_n$.]{
\label{f-KnKn:a}
\begin{minipage}[t]{0.99\textwidth}
\centering
\psfrag{p0}[bc]{$p_0$}
\psfrag{p1}[bl]{$p_1$}
\psfrag{p3}[bc]{$p_2$}
\psfrag{x}[tl]{$x$}
\psfrag{y}[tr]{$y$}
\psfrag{z}[bl]{$z$}
\psfrag{u}[tl]{$u$}
\psfrag{v}[tr]{$v$}
\includegraphics[width=0.55\textwidth]{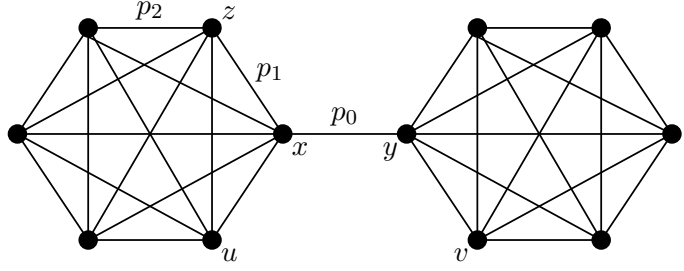}
\vspace{5ex}
\end{minipage}}
\vspace{7ex}

\subfigure[\normalsize Orbit chain under 
$C_2\ltimes(S_{n-1}\times S_{n-1})$.]{
\label{f-KnKn:b}
\begin{minipage}[t]{0.99\textwidth}
\centering
\psfrag{n-1p1}[bc]{$(n-1)p_1$}
\psfrag{p1}[tc]{$p_1$}
\psfrag{p0}[bc]{$p_0$}
\psfrag{Ox}[cl]{$O_x$}
\psfrag{Oz}[cr]{$O_z$}
\psfrag{x}[tl]{$x$}
\psfrag{y}[tr]{$y$}
\psfrag{Ou}[cr]{$O_u$}
\psfrag{Ov}[cl]{$O_v$}
\includegraphics[width=0.17\textwidth]{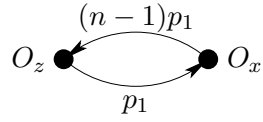}
\vspace{5ex}
\end{minipage}}
\vspace{7ex}

\subfigure[\normalsize Orbit chain under $S_{n-1}\times S_{n-1}$.]{
\label{f-KnKn:c}
\begin{minipage}[t]{0.99\textwidth}
\centering
\psfrag{n-1p1}[bc]{$(n-1)p_1$}
\psfrag{p1}[tc]{$p_1$}
\psfrag{p0}[bc]{$p_0$}
\psfrag{Ox}[cl]{$O_x$}
\psfrag{Oz}[cr]{$O_z$}
\psfrag{x}[tl]{$x$}
\psfrag{y}[tr]{$y$}
\psfrag{Ou}[cr]{$O_u$}
\psfrag{Ov}[cl]{$O_v$}
\includegraphics[width=0.4\textwidth]{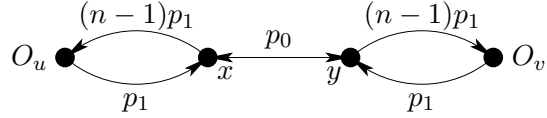}
\vspace{5ex}
\end{minipage}}
\vspace{7ex}

\subfigure[\normalsize Orbit chain under $S_{n-2}\times S_{n-1}$.]{
\label{f-KnKn:d}
\begin{minipage}[t]{0.99\textwidth}
\centering
\psfrag{p0}[bc]{$p_0$}
\psfrag{p1}[bl]{$p_1$}
\psfrag{np1}[tc]{$(n\!-\!2)p_1$}
\psfrag{np11}[tc]{$(n\!-\!1)p_1$}
\psfrag{p2}[tl]{$p_2$}
\psfrag{np2}[br]{$(n\!-\!2)p_2$}
\psfrag{x}[tl]{$x$}
\psfrag{y}[tr]{$y$}
\psfrag{z}[bl]{$z$}
\psfrag{Ou}[cr]{$O_u$}
\psfrag{Ov}[cl]{$O_v$}
\includegraphics[width=0.45\textwidth]{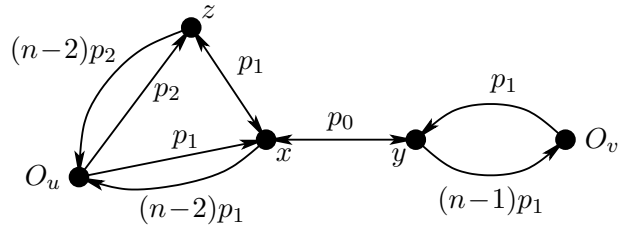}
\vspace{5ex}
\end{minipage}}
\vspace{7ex}

\caption{The graph $K_n$-$K_n$ and its orbit chains under 
different automorphism groups.
Here $O_x,O_z,O_u,O_v$ represent orbits of the vertices
$x,z,u,v$ (labeled in Figure~\ref{f-KnKn:a}), respectively, 
under the corresponding automorphism groups in each subgraph.}
\label{f-KnKn}
\end{figure}

The graph $K_n$-$K_n$ consists of two copies of the complete 
graph $K_n$ joined by a bridge (see Figure~\ref{f-KnKn:a}).
We follow the three steps described in~\S\ref{s-fmmc-reversible}.

First, it is clear by inspection that the full automorphism group
is $C_2\ltimes(S_{n-1}\times S_{n-1})$. 
The actions of $S_{n-1}\times S_{n-1}$ are all possible
permutations of the two set of $n-1$ vertices, distinct from the two
center vertices~$x$ and~$y$, among themselves.
The group $C_2$ acts on the graph by switching the two halves. 
The semi-direct product symbol $\ltimes$ means that the actions of 
$S_{n-1}\times S_{n-1}$ and $C_2$ do not commute. 

By symmetry analysis in~\S\ref{s-analysis}, there are three
edge orbits under the full automorphism group: the bridging edge
between vertices~$x$ and~$y$, the edges connecting~$x$ and~$y$ to all
other vertices, and the edges connecting all other vertices. 
Thus it suffices to consider just three transition probabilities
$p_0$, $p_1$, and $p_2$, each labeled in Figure~\ref{f-KnKn:a} 
on one representative of the three edge orbits.

As the second step, we construct the orbit chains.
The orbit chain of $K_n$-$K_n$ under the full automorphism group is
depicted in Figure~\ref{f-KnKn:b}. 
The orbit $O_x$ includes vertices~$x$ and~$y$, and the orbit $O_z$
consists of all other $2(n-1)$ vertices.
The transition probabilities of this orbit chain are calculated from
equation~(\ref{e-orbit-transition}) and are labeled on the directed
edges in Figure~\ref{f-KnKn:b}.
Similarly, the orbit chain under the subgroup $S_{n-1}\times S_{n-1}$
is depicted in Figure~\ref{f-KnKn:c}. 
While these two orbit chains are the most obvious to construct, 
none of them contains all eigenvalues of the original chain, nor
does their combination. 
For the one in Figure~\ref{f-KnKn:b}, the full automorphism
group does not have a fixed point either of its orbit $O_x$ or $O_z$. 
For the one in~\ref{f-KnKn:c}, the automorphism group
$S_{n-1}\times S_{n-1}$ has a fixed point in $O_x$ (either $x$ or
$y$), but does not have a fixed point in $O_z$ 
(note here $O_z$ is the orbit of~$z$ under the full automorphism
group). 
To fix the problem, we consider the orbit chain under the group 
$S_{n-2}\times S_{n-1}$, which leave the vertex~$x$, $y$, and $z$
fixed, while permuting the rest $n-2$ vertices on the left
and the $n-1$ points on the right, respectively. 
The corresponding orbit chain is shown in
Figure~\ref{f-KnKn:d}. 
By Corollary~\ref{c-full-eig-proj}, all distinct eigenvalues of the
original Markov chain on $K_n$-$K_n$ appear as eigenvalues of this
orbit chain.
Thus there are at most five distinct eigenvalues in the
original chain no matter how large~$n$ is.

To finish the second step, we calculate the transition probabilities
of the orbit chain under $H=S_{n-2}\times S_{n-1}$ using
equation~(\ref{e-orbit-transition}) and label them in Figure~\ref{f-KnKn:d}.
If we order the vertices of this orbit chain as 
$(x,y,z,O_u,O_v)$, then the transition probability
matrix on the orbit chain is 
\[
P_H=\left[ \begin{array}{ccccc}
1-p_0-(n-1)p_1 & p_0 & p_1 & (n-2)p_1 & 0 \\
p_0 & 1-p_0-(n-1)p_1 & 0 & 0 & (n-1)p_1 \\
p_1 & 0 & 1-p_1-(n-2)p_2 & (n-2)p_2 & 0 \\
p_1 & 0 & p_2 & 1-p_1-p_2 & 0 \\
0 & p_1 & 0 & 0 & 1-p_1 \end{array} \right] .
\]
By equation~(\ref{e-orbit-distribution}), the stationary
distribution of the orbit chain is 
\[
\pi_H=\left(\frac{1}{2n},~\frac{1}{2n},~\frac{1}{2n},
~\frac{n-2}{2n},~\frac{n-1}{2n} \right) .
\]

As the third step, we solve the SDP~(\ref{e-fmmc-rev})
with the above parametrization.
It is remarkable to see that we only need to solve an SDP with~$4$ 
variables (three transition probabilities $p_0$, $p_1$, $p_2$, 
and the extra scalar~$s$) and $5\times 5$ matrices 
no matter how large the graph (the number~$n$) is. 

We will revisit this example in \S\ref{s-Kn-Kn-2} using the block 
diagonalization method, where we present an analytic
expression for the exact optimal SLEM and corresponding transition
probabilities.

\subsubsection{Complete binary tree}
\label{s-tree-orbit}

\begin{figure}[t]
\centering
\subfigure[\normalsize Orbit graph and chain under 
$S_2 \wr S_2 \wr S_2$.]{
\label{f-binary-tree:a}
\begin{minipage}[t]{0.49\textwidth}
\centering
\psfrag{2p1}[br]{$2p_1$}
\psfrag{p1}[tl]{$p_1$}
\psfrag{2p2}[br]{$2p_2$}
\psfrag{p2}[tl]{$p_2$}
\psfrag{2p3}[br]{$2p_3$}
\psfrag{p3}[tl]{$p_3$}
\includegraphics[width=0.6\textwidth]{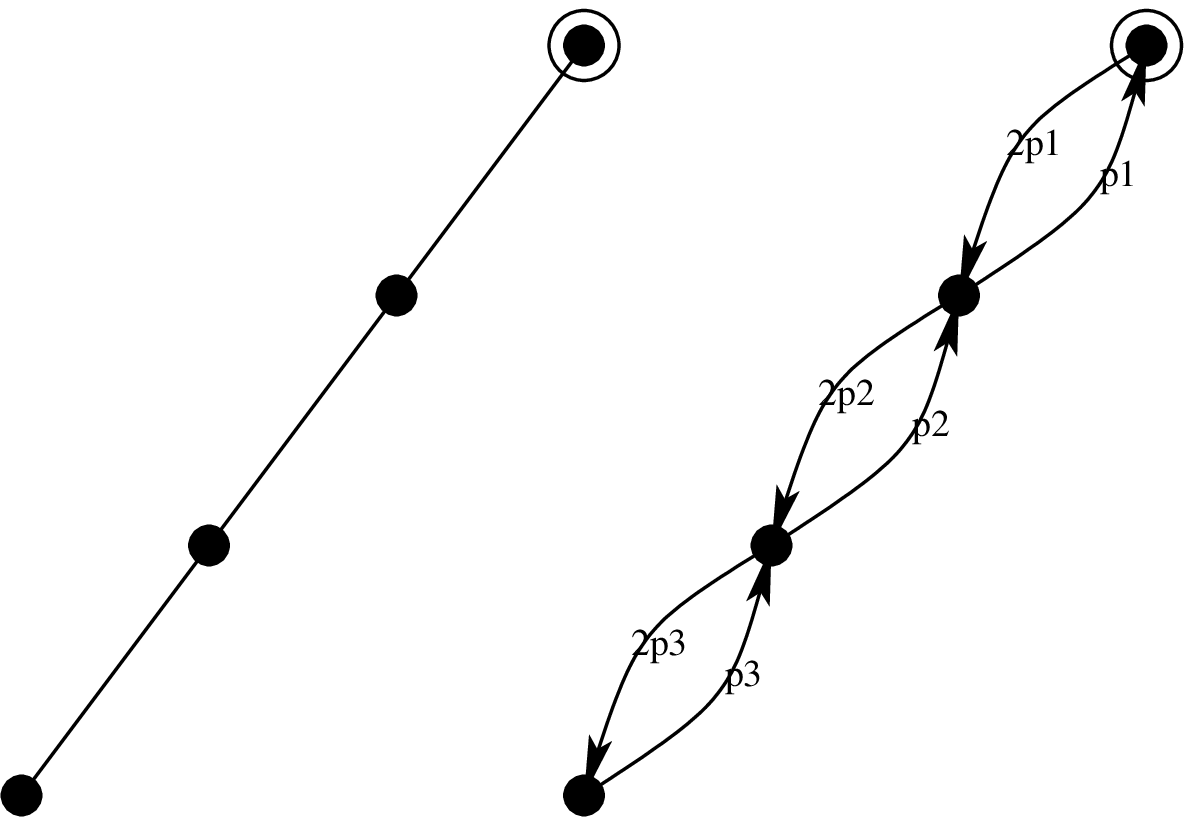}
\end{minipage}}%
\subfigure[\normalsize Orbit graph under 
$(S_2 \wr S_2) \times (S_2 \wr S_2)$.]{
\label{f-binary-tree:b}
\begin{minipage}[t]{0.49\textwidth}
\centering
\includegraphics[width=0.6\textwidth]{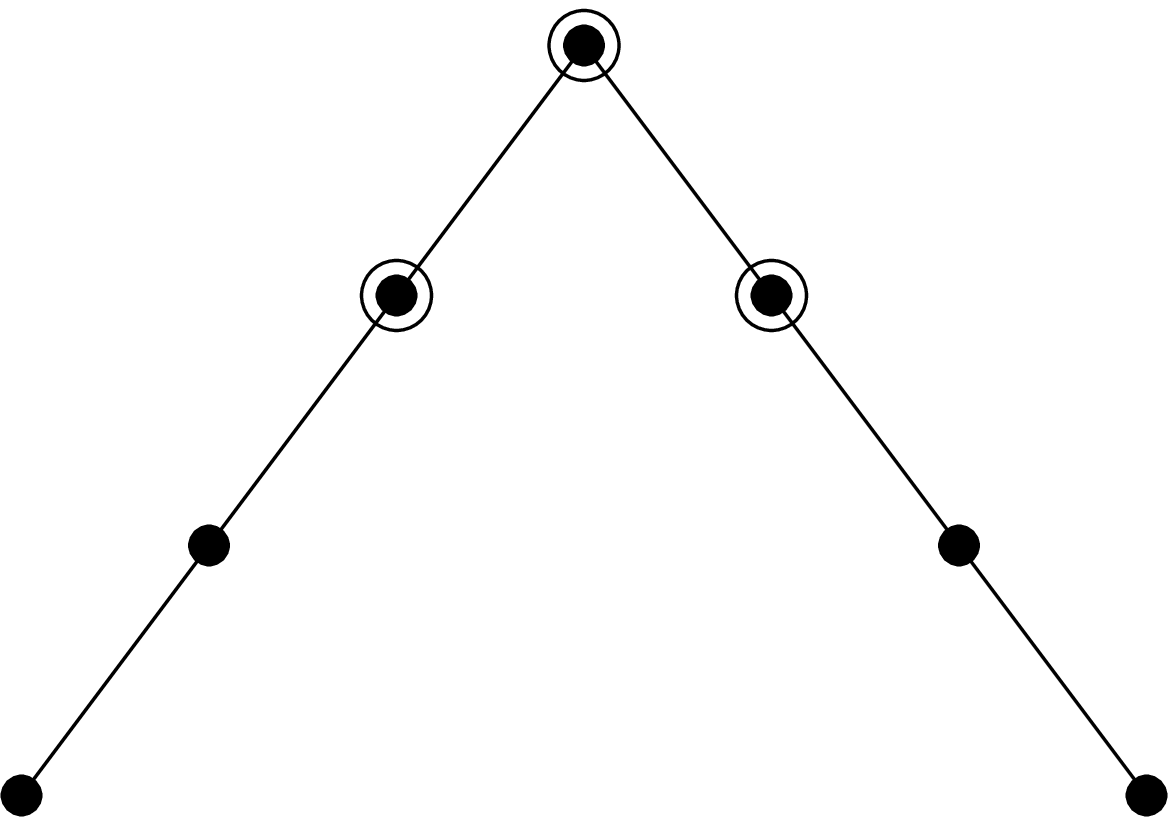}
\end{minipage}}
\newline\\[1ex]
\subfigure[\normalsize Orbit graph under 
$(S_2 \times S_2) \times (S_2 \wr S_2)$.]{
\label{f-binary-tree:c}
\begin{minipage}[t]{0.49\textwidth}
\centering
\includegraphics[width=0.6\textwidth]{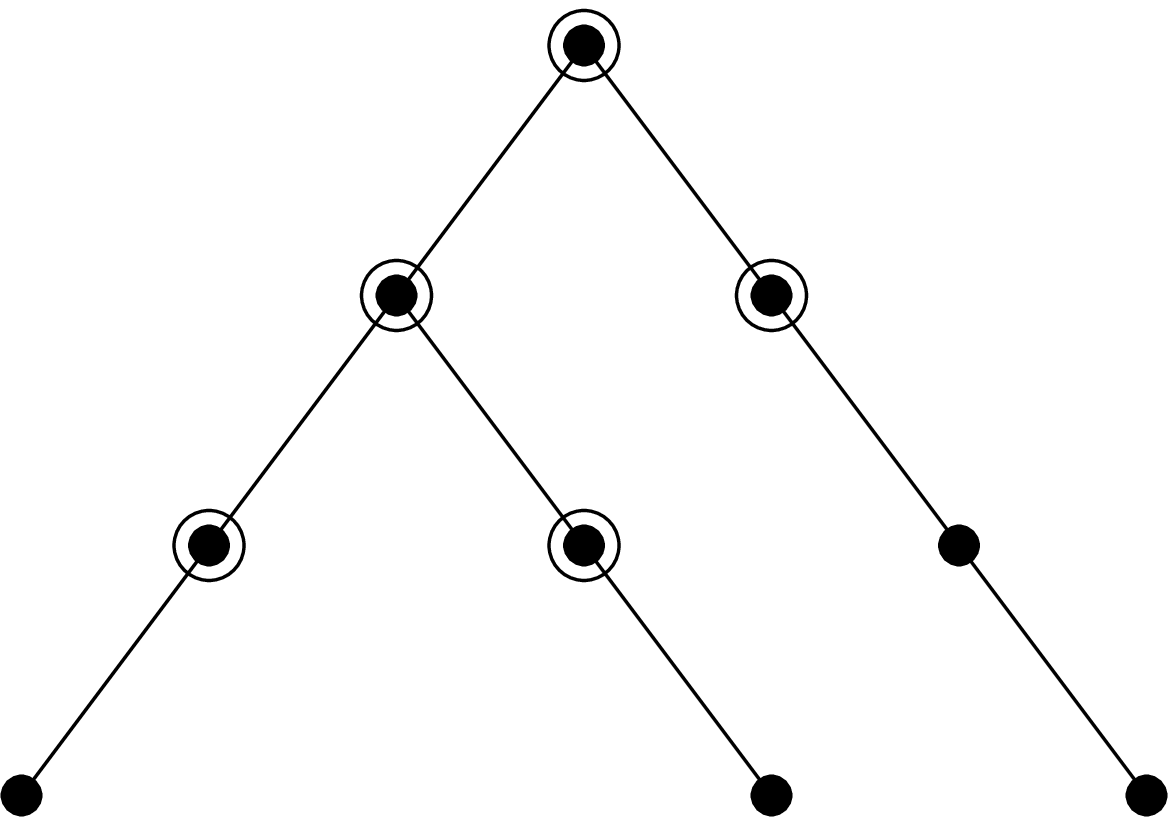}
\end{minipage}}%
\subfigure[\normalsize Orbit graph under 
$S_2 \times (S_2 \wr S_2)$.]{
\label{f-binary-tree:d}
\begin{minipage}[t]{0.49\textwidth}
\centering
\includegraphics[width=0.6\textwidth]{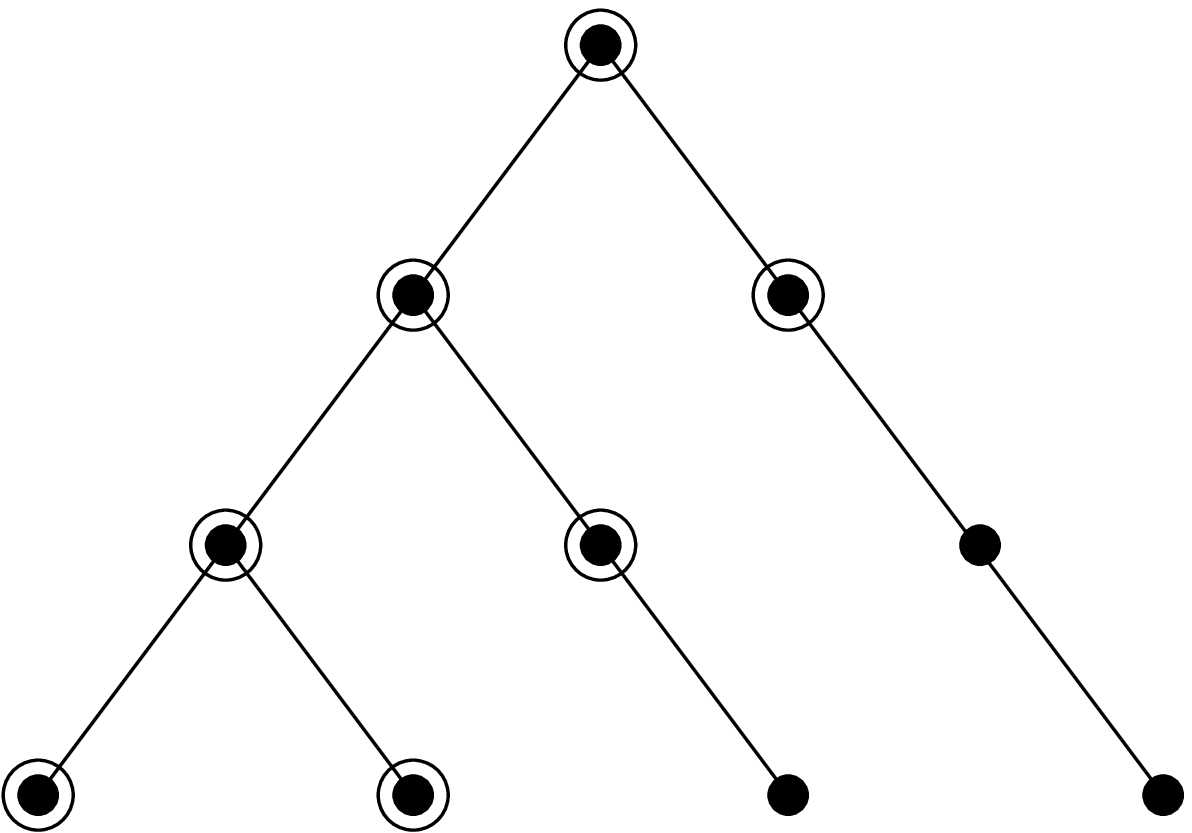}
\end{minipage}} 
\caption{Orbit graphs of the complete binary tree $\mathcal T_n$ 
($n=3$) under different automorphism groups. 
  The vertices surrounded by a circle are fixed
  points of the corresponding automorphism group.}
\label{f-binary-tree}
\end{figure}

We consider a complete binary tree with~$n$ levels of branches,
denoted as $\mathcal T_n$.  The total number of nodes is
$|\mathcal{V}|= 2^{n+1}-1$.  The matrix inequalities in the
corresponding SDP have size $|\mathcal{V}|\times|\mathcal{V}|$, which
is clearly exponential in $n$.  However, the binary tree has a very
large automorphism group, of size $2^{(2^n-1)}$.  This automorphism
group is best described recursively.  Plainly, for $n=1$, we have
$\mbox{Aut}(\mathcal T_1) = S_2$.  For $n>1$, it can be obtained by
the recursion
\[
\mbox{Aut}(\mathcal T_{k+1})=\mbox{Aut}(\mathcal T_k) \wr S_2, 
\qquad k=1,\ldots,n-1,
\]
where $\wr$ represents the \emph{wreath product} of two groups
(\eg, \cite{JaK:81}).
More specifically, let
$g=(g_1,g_2)$ and $h=(h_1,h_2)$ be elements of the product group 
$\mbox{Aut}(\mathcal T_k)\times\mbox{Aut}(\mathcal T_k)$, 
and $\sigma$ and $\pi$ be in $S_2$.
The multiplication rule of the wreath product is
\[
(g,\sigma)(h,\pi) = \left( (g_1 h_{\sigma^{-1}(1)}, 
g_2 h_{\sigma^{-1}(2)}), \sigma\pi \right).
\]
This is a semi-direct product $\mbox{Aut}(\mathcal T_k)^2\rtimes S_2$
(\cf~the automorphism group of $K_n$-$K_n$).
From the above recursion, the automorphism group of $\mathcal T_n$ is
\[
\mbox{Aut}(\mathcal T_n) = S_2 \wr S_2 \wr \cdots \wr S_2 
\quad \mbox{($n$ times)}.
\]
(The wreath product is associative, but not commutative.)
The representation theory of the automorphism group of the binary tree has been thoroughly studied as this group is the Sylow 2-subgroup of a symmetric group; see \cite{OOR:04,AbV:05}. 

The orbit graph of $\mathcal T_n$ under its full automorphism group is
a path with $n+1$ nodes (Figure~\ref{f-binary-tree:a}, left).
Since there are $n$ orbits of edges, there are $n$ different
transition probabilities we need to consider. 
We label them as $p_k$, $k=1,\ldots,n$, from top to bottom of the
tree. 
The corresponding orbit chain, represented by a directed graph labeled
with transition probabilities between orbits, is shown on the right of 
Figure~\ref{f-binary-tree:a}.
To simplify presentation, only the orbit graphs are shown in other
subfigures of Figure~\ref{f-binary-tree}.
The corresponding orbit chains should be straightforward to
construct. 

The largest subgroup of $\mbox{Aut}(\mathcal T_n)$ that has a fixed
point in every orbit under $\mbox{Aut}(\mathcal T_n)$ is
\[
W_n = \prod_{k=1}^{n-1} (S_2 \wr \cdots \wr S_2)~(k~\mbox{times})
\]
where $\prod$ denotes direct product of groups.  
The corresponding orbit graph is shown in Figure~\ref{f-binary-tree:d}
for $n=3$. 
The number of vertices in this orbit graph is 
\[
1+2+\cdots+n+(n+1) = {n+1 \choose 2} = \frac{1}{2}(n+1)(n+2),
\] 
which is much smaller than $2^{n+1}-1$, the size of $\mathcal T_n$.

From the above analysis, we only need to solve the fastest reversible
Markov chain problem on the orbit graph of size ${n+1 \choose 2}$ with
$n$ variables $p_1,\ldots,p_n$.  In next section, using the technique
of block diagonalization, we will see that the transition probability
matrix of size ${n+1 \choose 2}$ can be further decomposed into smaller
matrices with sizes $1,2,\ldots,n+1$.  Due to an eigenvalue
interlacing result, we only need to consider the orbit chain with
$2n+1$ vertices in Figure~\ref{f-binary-tree:b}.  

\section{Symmetry reduction by block diagonalization}
\label{s-block-diag}

By definition of the fixed-point subspace $\mathcal F$
(in~\S\ref{s-fmmc-invariant}), any transition probability matrix
$P\in\mathcal F$ is invariant under the actions of
$\mbox{Aut}(\mathcal G)$. 
More specifically, for any permutation matrix $Q$ given by
$\sigma\in\mbox{Aut}(\mathcal G)$, we have $Q P Q^T=P$, equivalently 
$Q P = P Q$.
In this section we show that this property allows the construction 
of a coordinate transformation matrix that can block diagonalize 
every $P\in\mathcal F$. 
The resulting blocks usually have much smaller sizes and repeated
blocks can be discarded in computation.

The method we use in this section is based on classical group
representation theory (\eg, \cite{Ser:77}).
It was developed for more general SDPs in~\cite{GatermannParrilo}, 
and has found applications in sum-of-squares decomposition for
minimizing polynomial
functions~\cite{PhD:Parrilo,Par:03,ParriloSturmfels} and controller
design for symmetric dynamical systems~\cite{CLP:03}.
A closely related approach is developed in \cite{dPS:07}, which is 
based on a low-order representation of the commutant 
(collection of invariant matrices) of the matrix algebra generated 
by the permutation matrices. 

\subsection{Some group representation theory}
\label{s-representation}
Let $G$ be a group.
A \emph{representation} $\rho$ of $G$ assigns an invertible matrix
$\rho(g)$ to each $g\in G$ in such a way that the matrix assigned to
the product of two elements in $G$ is the product of the matrices
assigned to each element: $\rho(gh)=\rho(g)\rho(h)$.
The matrices we work with are all invertible and are considered over
the real or complex numbers.
We thus regard $\rho$ as a homomorphism from $g$ to the linear maps on
a vector space $V$.
The dimension of $\rho$ is the dimension of $V$.
Two representations are \emph{equivalent} if they are related by a
fixed similarity transformation. 

If $W$ is a subspace of $V$ invariant under $G$, then $\rho$
restricted to $W$ gives a \emph{subrepresentation}.
Of course the zero subspace and the subspace $W=V$ are trivial
subrepresentations. 
If the representation $\rho$ admits no non-trivial subrepresentation,
then $\rho$ is called \emph{irreducible}.

We consider first complex representations, as the theory is
considerably simpler in this case. 
For a finite group $G$ there are only finitely many inequivalent
irreducible representations $\vartheta_1,\ldots,\vartheta_h$ of 
dimensions $n_1,\ldots,n_h$, respectively. 
The degrees $n_i$ divide the group order $|G|$, and satisfy the
condition $\sum_{i=1}^h n_i^2=|G|$.
Every linear representation of $G$ has a canonical decomposition as a
direct sum of irreducible representations 
\[
\rho = m_1\vartheta_1 \oplus m_2\vartheta_2 \oplus \cdots \oplus
m_h\vartheta_h, 
\]
where $m_1,\ldots,m_h$ are the multiplicities.
Accordingly, the representation space $\mathbf{C}^n$ has an
\emph{isotypic decomposition} 
\begin{equation}\label{e-isotypic}
\mathbf{C}^n = V_1 \oplus \cdots \oplus V_h
\end{equation}
where each isotypic components consists of $m_i$ invariant subspaces
\begin{equation}\label{e-equivalent}
V_i = V_i^1 \oplus \cdots \oplus V_i^{m_i},
\end{equation}
each of which has dimension $n_i$ and transforms after the manner of
$\vartheta_i$. 
A basis of this decomposition transforming with respect to the
matrices $\vartheta_i(g)$ is called \emph{symmetry-adapted} and can be
computed using the algorithm presented in \cite[\S2.6-2.7]{Ser:77} 
or \cite[\S5.2]{FaS:92}. 
This basis defines a change of coordinates by a matrix $T$ collecting
the basis as columns.
By Schur's lemma, if a matrix~$P$ satisfies 
\begin{equation}\label{e-invariance}
\rho(g)P=P\rho(g), \qquad \forall g\in G,
\end{equation}
then $T^{-1}PT$ has block diagonal form with one block
$P_i$ for each isotypic component of dimension $m_i n_i$, which
further decomposes into $n_i$ equal blocks $B_i$ of dimension $m_i$. 
That is
\begin{equation}\label{e-block-diag}
T^{-1} P T = \left[ \begin{array}{ccc}
P_1 & & 0\\ & \ddots & \\ 0 & & P_h \end{array} \right], \qquad
P_i = \left[ \begin{array}{ccc}
B_i & & 0\\ & \ddots & \\ 0 & & B_i \end{array} \right].
\end{equation}

For our application of semidefinite programs, the problems are
presented in terms of real matrices, and therefore we would like to
use real coordinate transformations. 
In fact a generalization of the classical theory to the real case is
presented in \cite[\S13.2]{Ser:77}.
If all $\vartheta_i(g)$ are real matrices the irreducible
representation is called \emph{absolutely irreducible}. 
Otherwise, for each $\vartheta_i$ with complex character its complex
conjugate will also appear in the canonical decomposition.
Since $\rho$ is real both will have the same multiplicity and real
bases of $V_i+\bar V_i$ can be constructed. 
So two complex conjugate irreducible representations form one real
irreducible representation of \emph{complex type}. 
There is a third case, real irreducible representations of
\emph{quaternonian type}, rarely seen in practical examples. 

In this paper, we assume that the representation $\rho$ is orthogonal,
\ie, $\rho(g)^T\rho(g)=\rho(g)\rho(g)^T=I$ for all $g\in G$.
As a result, the transformation matrix~$T$ can also be chosen to be
orthogonal.
Thus $T^{-1}=T^T$ (for complex matrices, it is the conjugate
transpose).  
For symmetric matrices the block corresponding to a representation of
complex type or quaternonian type simplifies to a collection of equal
subblocks. 
For the special case of circulant matrices, complete
diagonalization reveals all the eigenvalues \cite[page 50]{Dia:88}.

\subsection{Block diagonalization of SDP constraint}
\label{s-diag-sdp}

As in~\S\ref{s-fmmc-invariant}, for every 
$\sigma\in\mbox{Aut}(\mathcal G)$ we assign a permutation matrix
$Q(\sigma)$ by letting $Q_{ij}(\sigma)=1$ if 
$i=\sigma(j)$ and $Q_{ij}(\sigma)=0$ otherwise.  
This is an $n$-dimensional representation of $\mbox{Aut}(\mathcal G)$,
which is often called the \emph{natural representation}.
As mentioned in the beginning of this section, every matrix~$P$ 
in the fixed-point subset $\mathcal F$ has the symmetry of 
$\mbox{Aut}(\mathcal G)$; \ie, it satisfies the
condition~(\ref{e-invariance}) with $\rho=Q$. 
Thus a coordinate transformation matrix~$T$ can be constructed such
that $P$ can be block diagonalized into the
form~(\ref{e-block-diag}). 

Now we consider the SDP~(\ref{e-fmmc-orbits-sdp}), which is the 
FMMC problem formulated in the fixed-point subset~$\mathcal F$.
In \S\ref{s-reduce-variables}, we have derived the expression 
$P(p)=I-\sum_{k=1}^N p_k L_k$,
where $L_k$ is the Laplacian matrix for the $k$th orbit graph and $p_k$ is the common transition probability assigned on all edges in the $k$th orbit graph. 
Note the matrix $P(p)$ has the symmetry of $\mbox{Aut}(\mathcal G)$.
Applying the coordinate transformation~$T$ to the linear matrix
inequalities, we obtain the following equivalent problem
\begin{equation}\label{e-fmmc-block-diag}
\begin{array}{ll}
\mbox{minimize}   & s \\[1ex]
\mbox{subject to} & -s I_{m_i} \preceq B_i(p) - J_i \preceq s I_{m_i},
\quad i=1,\ldots,h\\[1ex]
& p_k \geq 0, \quad k=1,\ldots,N\\[1ex]
& \sum_{k=1}^N (L_k)_{ii} \;p_k \leq 1, \quad i=1,\ldots,n
\end{array}
\end{equation}
where $B_i(p)$ correspond to the small blocks $B_i$
in~(\ref{e-block-diag}) of the transformed matrix $T^TP(p)T$, 
and~$J_i$ are the corresponding diagonal blocks of 
$T^T(1/n)\ones\ones^T T$.
The number of matrix inequalities~$h$ is the number of inequivalent
irreducible representations, and the size of each matrix
inequality~$m_i$ is the multiplicity of the corresponding 
irreducible representation. 
Note that we only need one out of~$n_i$ copies of each~$B_i$ in the
decomposition~(\ref{e-block-diag}).
Since $m_i$ can be much smaller than~$n$ (the number of vertices in
the graph), the improvement in computational complexity over the SDP
formulation~(\ref{e-fmmc-orbits-sdp}) can be significant (see the flop 
counts discussed in~\S\ref{s-structure}).    
This is especially the case when there are high-dimensional
irreducible representations 
(\ie, when $n_i$ is large; see, \eg, $K_n$-$K_n$ defined in~\S\ref{s-Kn-Kn}). 

The transformed SDP formulation~(\ref{e-fmmc-block-diag}) needs some
further justification. 
Namely, all the off-diagonal blocks of the matrix
$T^T(1/n)\ones\ones^T T$ have to be zero.
This is in fact the case.
Moreover, the following theorem reveals an interesting connection
between the block diagonalization approach and the orbit theory
in~\S\ref{s-fmmc-orbit}. 

\begin{theorem}\label{t-orbit-diag-connection}
Let $H$ be a subgroup of $\mathrm{Aut}(\mathcal G)$, and $T$ be the
coordinate transformation matrix whose columns are a symmetry-adapted
basis for the natural representation of~$H$.  Suppose a Markov
chain~$P$ defined on the graph has the symmetry of~$H$.  Then the
matrix $T^T(1/n)\ones\ones^T T$ has the same block diagonal form as
$T^TPT$.  Moreover, there is only one nonzero block.  Without loss of
generality, let this nonzero block be $J_1$ and the corresponding
block of $T^TPT$ be $B_1$.  These two blocks relate to the orbit chain
$P_H$ by
\begin{eqnarray}
B_1 &=& \Pi^{1/2}P_H\Pi^{-1/2}
        \label{e-orbit-diag-connection} \\
J_1 &=& qq^T \label{e-orbit-diag-J1}
\end{eqnarray}
where $\Pi=\Diag(\pi_H)$, $q=\sqrt{\pi_H}$,
and $\pi_H$ is the stationary distribution of $P_H$.  
\end{theorem}

\begin{proof}
First we note that $P$ always has a single eigenvalue~$1$ with
associated eigenvector~$\ones$.  Thus~$\ones$ spans an invariant
subspace of the natural representation, which is obviously
irreducible.  The corresponding irreducible representation is
isomorphic to the trivial representation (which assigns the scalar~$1$
to every element in the group).  Without loss of generality, let $V_1$
be the isotypic component that contains the vector~$\ones$.  Thus
$V_1$ is a direct product of $H$-fixed vectors (each corresponds to a
copy of the trivial representation), and~$\ones$ is a linear
combination of these vectors.

Let $m_1$ be the dimension of $V_1$, which is the number of $H$-fixed
vectors.  We can calculate~$m_1$ by Frobenius reciprocity, or
``Burnside's Lemma''; see, \eg, \cite{Ser:77}.  To do so, we note that
the character~$\chi$ of the natural representation~$Q(g)$,
$g\in H$, is the number of fixed points of~$g$, \ie,
\[
\chi(g) = \Tr Q(g) = \mathrm{FP}(g)=\#\{v\in \mathcal V: g(v)=v \}.
\]
``Burnside's Lemma'' says that
\[
\frac{1}{|H|} \sum_{g\in H} \mathrm{FP}(g) = \# \mathrm{orbits}. 
\]
The left-hand side is the inner product of~$\chi$ with the trivial
representation.  It thus counts the number of $H$-fixed vectors
in~$V$.  So $m_1$ equals the number of orbits under $H$.

Suppose that $\mathcal V=O_1\cup\ldots\cup O_{m_1}$ as a disjoint
union of $H$-orbits.
Let $b_i(v)=1/\sqrt{|O_i|}$ if $v\in O_i$ and zero otherwise.  
Then $b_1,\ldots,b_{m_1}$ are $H$-fixed
vectors, and they form an orthonormal symmetry-adapted basis 
for $V_1$ (these are not unique).  
Let $T_1=[b_1 \cdots b_{m_1}]$ be the first $m_1$ columns 
of~$T$. 
They are orthogonal to all other columns of~$T$.  
Since~$\ones$ is a linear combination of $b_1,\ldots,b_{m_1}$, it is
also orthogonal to other columns of~$T$. 
Therefore the matrix $T^T(1/n)\ones\ones^T T$ has all its elements
zero except for the first $m_1\times m_1$ diagonal block, which we
denote as $J_1$.
More specifically, $J_1=qq^T$ where 
\begin{eqnarray*}
q &=&\frac{1}{\sqrt{n}} T_1^T \ones = \frac{1}{\sqrt{n}} 
     \left[ b_1^T\ones ~\cdots~ b_{m_1}^T\ones \right]^T\\
&=& \frac{1}{\sqrt{n}} \left[ \frac{|O_1|}{\sqrt{|O_1|}} ~\cdots~  
    \frac{|O_{m_1}|}{\sqrt{|O_{m_1}|}} \; \right]^T
    = \left[ \sqrt{\frac{|O_1|}{n}} ~\ldots~ 
      \sqrt{\frac{|O_{m_1}|}{n}} \; \right]^T.
\end{eqnarray*}
Note that by~(\ref{e-orbit-distribution}) the stationary distribution 
of the orbit chain $P_{H}$ is 
\[
\pi_{H}
=\left[ \; \frac{|O_1|}{n} ~\cdots~ \frac{|O_{m_1}|}{n} \; \right]^T.
\]   
Thus we have $q=\sqrt{\pi_{H}}$.
This proves~(\ref{e-orbit-diag-J1}).

Finally we consider the relationship between $B_1=T_1^T P T_1$ and 
$P_{H}$. 
We prove~(\ref{e-orbit-diag-connection}) by showing
\[
\Pi^{-1/2} B_1 \Pi^{1/2} = \Pi^{-1/2} T_1^T P T_1 \Pi^{1/2} = P_H. 
\] 
It is straightforward to verify that
\begin{eqnarray*}
&&
\Pi^{-1/2} T_1^T = \sqrt{n} \left[ \begin{array}{c}
b'^T_1\\ \vdots\\ b'^T_{m_1} \end{array} \right], \qquad
b'_i(v) = \left\{ \begin{array}{lll} 
          \displaystyle \frac{1}{|O_i|} && \mbox{if}~v\in O_i\\[2ex]
          0 && \mbox{if}~v\notin O_i \end{array} \right. \\[2ex]
&&
T_1 \Pi^{1/2} = \frac{1}{\sqrt{n}} 
                \left[ b''_1 ~\cdots~ b''_{m_1}\right], \qquad
b''_i(v) = \left\{ \begin{array}{lll} 
           1 && \mbox{if}~v\in O_i\\[1ex]
           0 && \mbox{if}~v\notin O_i \end{array} \right.
\end{eqnarray*}
The entry at the $i$-th row and $j$-th column of the matrix 
$\Pi^{-1/2} T_1^T P T_1 \Pi^{1/2}$ are given by
\[
b'^T_i P b''_j 
= \frac{1}{|O_i|} \sum_{v\in O_i} \sum_{u\in O_j} P(v,u)
= \frac{1}{|O_i|} \sum_{v\in O_i} P_H(O_i, O_j)
= P_H( O_i, O_j ).
\]
In the last equation, we have used the fact that $P_H(O_i, O_j)$ is
independent of which $v\in O_i$ is chosen.
This completes the proof.
\end{proof}

From Theorem~\ref{t-orbit-diag-connection}, we know that
$B_1$ contains the eigenvalues of the orbit chain under $H$. 
Other blocks $B_i$ contain additional eigenvalues (not including those
of $P_{H}$) of the orbit 
chains under various subgroups of $H$.
(Note that the eigenvalues of the orbit chain under
$H$ are always contained in the orbit chain 
under its subgroups).
With this observation, it is possible to identify the multiplicities
of eigenvalues in orbit chains under various subgroups of 
$\mbox{Aut}(\mathcal{G})$ by relating to the
decompositions~(\ref{e-isotypic}), (\ref{e-equivalent}) and
(\ref{e-block-diag})
(some preliminary results are discussed in \cite{BDPX:05}). 

\subsubsection{A running example}

\begin{figure}
\begin{center}
\psfrag{1}{$1$}
\psfrag{2}{$2$}
\psfrag{3}{$3$}
\psfrag{4}{$4$}
\psfrag{5}{$5$}
\psfrag{6}{$6$}
\psfrag{7}{$7$}
\psfrag{8}{$8$}
\psfrag{9}{$9$}
\includegraphics[width=0.2\textwidth]{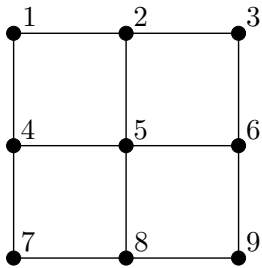}
\caption{A $3 \times 3$ grid graph.}
\label{fig:grid9}
\end{center}
\end{figure}

As a running example for this section, we consider a
Markov chain on a $3\times 3$ grid $\mathcal{G}$, with a total of 9
nodes (see Figure~\ref{fig:grid9}). The automorphism group
$\mbox{Aut}(\mathcal{G})$ is isomorphic to the 8-element dihedral group
$D_4$, and corresponds to flips and 90-degree rotations of the graph.
The orbits of $\mbox{Aut}(\mathcal{G})$ acting on the vertices and
edges are
\[
\{1,3,7,9\}, \qquad \{5\}, \qquad \{2,4,6,8\}
\]  
and 
\[
\{ \{1,2\}, \{1,4\}, \{2,3\}, \{3,6\}, \{4,7\}, \{7,8\}, 
   \{6,9\}, \{8,9\} \}, \qquad
\{ \{2,5\}, \{4,5\}, \{5,6\}, \{5,8\} \},
\]
respectively. 
So $\mathcal{G}$ is neither vertex- nor edge-transitive.

\begin{figure}
\begin{center}
\psfrag{o1}[cr]{$O_1$}
\psfrag{o2}[bc]{$O_2$}
\psfrag{o5}[cl]{$O_5$}
\psfrag{2aa}[bc]{$2a$}
\psfrag{2ab}[tc]{$2a$}
\psfrag{b}[bc]{$b$}
\psfrag{4b}[tc]{$4b$}
\includegraphics[width=0.35\textwidth]{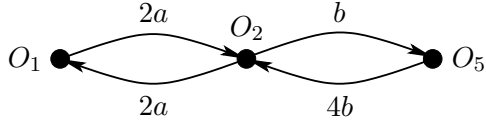}
\caption{The orbit chain of the $3 \times 3$ grid graph.}
\label{f-grid-orbit}
\end{center}
\end{figure}

By Corollary~\ref{c-numprobs}, we associate transition
probabilities~$a$ and~$b$ to the two edge orbits, respectively. 
The transition probability matrix has the form
\[
P=\left[\begin{array}{ccccccccc}
1\!-\!2a & a & 0 & a & 0 & 0 & 0 & 0 & 0 \\
a & 1\!-\!2a\!-\!b & a & 0 & b & 0 & 0 & 0 & 0 \\
0 & a & 1\!-\!2a & 0 & 0 & a & 0 & 0 & 0 \\
a & 0 & 0 & 1\!-\!2a\!-\!b & b & 0 & a & 0 & 0\\
0 & b & 0 & b & 1\!-\!4b & b & 0 & b & 0 \\
0 & 0 & a & 0 & b & 1\!-\!2a\!-\!b & 0 & 0 & a \\
0 & 0 & 0 & a & 0 & 0 & 1\!-\!2 a & a & 0 \\
0 & 0 & 0 & 0 & b & 0 & a & 1\!-\!2a\!-\!b& a \\
0 & 0 & 0 & 0 & 0 & a & 0 & a & 1\!-\!2a 
\end{array}
\right].
\]
The matrix $P$ satisfies $Q(\sigma)P=P Q(\sigma)$ for every
$\sigma\in\mbox{Aut}(\mathcal G)$.  Using the algorithm in
\cite[\S5.2]{FaS:92}, we found a symmetry-adapted basis for the
representation $Q$, which we take as columns to form
\[
T = \frac{1}{2}
\left[\begin{array}{rrrrrrrrr}
0 & 1 & 0 & 1 & 0 & \sqrt{2} & 0     & 0      & 0 \\
0 & 0 & 1 & 0 &-1 & 0      & 1 & 0      & 1 \\
0 & 1 & 0 &-1 & 0 & 0      & 0     & \sqrt{2} & 0 \\
0 & 0 & 1 & 0 & 1 & 0      & 1 & 0      & -1\\
2 & 0 & 0     & 0 & 0     & 0      & 0     & 0      & 0 \\
0 & 0 & 1 & 0 & 1 & 0      &-1 & 0      & 1 \\
0 & 1 & 0 &-1 & 0 & 0      & 0     & -\sqrt{2} & 0 \\
0 & 0 & 1 & 0 &-1 & 0      &-1 & 0      & -1 \\
0 & 1 & 0 & 1 & 0 & -\sqrt{2} & 0     & 0      & 0
\end{array}
\right].
\]
With this coordinate transformation matrix, we obtain
\[
T^TPT=\left[\begin{array}{ccccccccc}
1\!-\!4 b & 0 & 2b &  &   &   &   &   &   \\
0 & 1\!-\!2a & 2a &   &   &   &   &   &   \\
2b  & 2 a & 1\!-\!2a\!-\!b &   &   &   &   &   &   \\
  &   &   & 1\!-\!2a &   &   &   &   &  \\
  &   &   &   & 1\!-\!2a\!-\!b &   &   &   &   \\
  &   &   &   &   & 1\!-\!2a & \sqrt{2} a &   &   \\
  &   &   &   &   & \sqrt{2} a & 1\!-\!2 a \!-\!b &   &   \\
  &   &   &   &   &   &   & 1\!-\!2a & \sqrt{2} a \\
  &   &   &   &   &   &   & \sqrt{2} a & 1\!-\!2a\!-\!b 
\end{array}
\right].
\]
The 3-dimensional block $B_1$ contains the single eigenvalue~$1$, and
it is related to the orbit chain in Figure~\ref{f-grid-orbit} by the
equation~(\ref{e-orbit-diag-connection}). 
The corresponding nonzero block of $T^T(1/n)\ones\ones^T T$ is
\[
J_1 = \frac{1}{9} \left[ \begin{array}{ccc}
1 & 2 & 2\\ 2 & 4 & 4\\ 2 & 4 & 4 \end{array} \right].
\]

Next, we substitute the above expressions into the
SDP~(\ref{e-fmmc-block-diag}) and solve it numerically.
Since there are repeated $2\times 2$ blocks, the original $9 \times 9$
matrix is replaced by four smaller blocks, of dimension 3,1,1,2.
The optimal solutions are
\[
a^\star \approx 0.363, \qquad b^\star \approx 0.2111, \qquad 
\mu^\star \approx 0.6926.
\]
Interestingly, it can be shown that these optimal values are not
rational, but instead algebraic numbers with defining minimal
polynomials:
\begin{eqnarray*}
18157\, a^5-17020\, a^4+6060\, a^3-1200\, a^2+180\, a-16 &=& 0 \\
1252833\, b^5-1625651\, b^4+791936\, b^3-173536\, b^2+15360\, b-256 &=& 0 \\
54471\, \mu^5-121430\, \mu^4+88474\, \mu^3-18216\, \mu^2-2393\, \mu+262 &=& 0.
\end{eqnarray*}


\subsection{Examples}
We revisit some previous examples with the block diagonalization
method, and draw connections to the method based on orbit theory
in~\S\ref{s-fmmc-orbit}. 
We also discuss some additional examples that are difficult if one uses 
the orbit theory, but are nicely handled by block diagonalization.  
In many of the examples, the coordinate transformation matrix~$T$ 
can be constructed directly by inspection.

\subsubsection{Complete bipartite graphs}
\label{s-bipartite-diag}
For the complete bipartite graph $K_{m,n}$ (see
Figure~\ref{f-bipartite}),  
This graph is edge-transitive, so we can assign the same transition
probability~$p$ on all the edges. 
The transition probability matrix has the form
\[
P(p)=\left[
\begin{array}{cc}
(1 - n p) I_m & p \, \ones_{m\times n} \\
p \, \ones_{n\times m} & (1- m p) I_n 
\end{array}
\right]
\]
We can easily find a decomposition of the associated matrix
algebra. It will have three blocks, and an orthogonal
block-diagonalizing change of basis is given by
\[
T = \left[ \begin{array}{cccc}
(1/\sqrt{m})\ones_{m\times 1} & 0 & F_m & 0 \\
0 & (1/\sqrt{n}) \ones_{n\times 1} & 0 & F_m          
\end{array} \right]
\]
where $F_n$ is an $n \times (n-1)$ matrix whose columns are an
orthogonal basis of the subspace complementary to that generated by 
$\ones_{n \times 1}$.

In the new coordinates, the matrix $T^T P(p) T$ has the following
diagonal blocks
\[
\left[
\begin{array}{cc}
1-m p & p \sqrt{nm} \\
p \sqrt{nm} & 1-n p
\end{array}
\right],
\qquad
I_{n-1} \otimes (1-m p), \qquad
I_{m-1} \otimes (1-n p).
\]
The $2\times 2$ block has eigenvalues~$1$ and~$1-(m+n)p$. 
The other diagonals reveal the eigenvalue $1-mp$ and $1-np$, with
multiplicities $n-1$ and $m-1$, respectively. 
The optimal solution to the FMMC problem can be easily obtained as
in~(\ref{e-bipartite-alpha}) and~(\ref{e-bipartite-mu}).

To draw connections to the orbit theory, we note that the above
$2\times 2$ block is precisely $B_1$ in the 
equation~(\ref{e-orbit-diag-connection}), and the corresponding 
$P_H$ is the orbit chain shown in Figure~\ref{f-bipartite-orbit:a}.
In addition to the two eigenvalues in~$B_1$, the extra eigenvalue in
the orbit chain of Figure~\ref{f-bipartite-orbit:b} is $1-np$, and the 
extra eigenvalue in Figure~\ref{f-bipartite-orbit:c} is $1-mp$. 
All these eigenvalues appear in the orbit chain in 
Figure~\ref{f-bipartite-orbit:d}.
As we have seen, the block diagonalization technique reveals the
multiplicities in the original chain of the eigenvalues from various
orbit chains.

\subsubsection{Complete $k$-partite graphs}
The previous example generalizes nicely to the complete $k$-partite
graph $K_{n_1,\ldots,n_k}$. 
In this case, the fixed-point reduced matrix will have
dimensions $\sum _i n_i$, and the structure
\[
P(p)=\left[
\begin{array}{cccc}
(1 -  \sum_{j \not = 1} n_j p_{1j} ) I_{n_1} & 
p_{12} \ones_{n_1 \times n_2} & \cdots & 
p_{1k} \ones_{n_1 \times n_k}\\
p_{21} \ones_{n_2 \times n_1} & 
(1 -  \sum_{j \not = 2} n_j p_{2j} ) I_{n_2} & 
\cdots & 
p_{2k} \ones_{n_2 \times n_k}\\
\vdots & \vdots & \ddots & \vdots \\
p_{k1} \ones_{n_k \times n_1} & 
p_{k2} \ones_{n_k \times n_2} &
\cdots & 
(1 -  \sum_{j \not = k} n_j p_{kj} ) I_{n_k} 
\end{array}
\right]
\]
where the probabilities satisfy $p_{ij}=p_{ji}$.
There are total ${k \choose 2}$ independent variables.

In a very similar fashion to the bipartite case, we can explicitly
write the orthogonal coordinate transformation matrix
\[
T = \left[ \begin{array}{cccccc}
(1/\sqrt{n_1}) \ones_{n_1 \times 1}   & \ldots       & \mathbf{0} &
F_{n_1}    & \ldots     & \mathbf{0}  \\
\vdots & \ddots & \vdots & \vdots & \ddots & \vdots \\
\mathbf{0} &  \ldots & (1/\sqrt{n_{k}}) \ones_{n_{k} \times 1}  &
\mathbf{0} & \ldots  &   F_{n_k}          
\end{array} \right].
\]
The matrix $T^T P(p) T$ decomposes into $k+1$ blocks: one of
dimension $k$, with the remaining $k$ blocks each having dimension
$n_i-1$. The decomposition is:
\begin{eqnarray*}
&& \left[ \begin{array}{cccc}
(1 -  \sum_{j \not = 1} n_j p_{1j} ) & 
p_{12} \sqrt{n_1 n_2}& \cdots & 
p_{1k} \sqrt{n_1 n_k} \\
p_{21} \sqrt{n_2 n_1}& 
(1 -  \sum_{j \not = 2} n_j p_{2j} ) & 
\cdots & 
p_{2k} \sqrt{n_2 n_k} \\
\vdots & \vdots & \ddots & \vdots \\
p_{k1} \sqrt{n_k n_1} & 
p_{k2} \sqrt{n_k n_2} &
\cdots & 
(1 -  \sum_{j \not = k} n_j p_{kj} ) 
\end{array} \right], \\[2ex]
&& I_{n_i-1} \otimes (1 -\sum_{j \not = i} n_j p_{i j} ), 
\qquad i = 1,\ldots,k.
\end{eqnarray*}
These blocks can be substituted into the SDP~(\ref{e-fmmc-block-diag})
to solve the FMMC problem.

\subsubsection{Wheel graph}

\begin{figure}
\begin{center}
\includegraphics[width=0.22\textwidth]{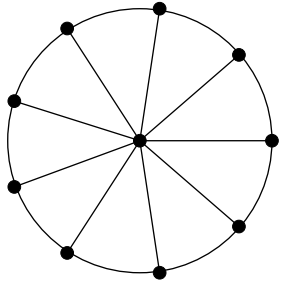}
\caption{The wheel graph with $n=9$ (total 10 nodes).}
\label{f-wheel}
\end{center}
\end{figure}

The wheel graph consists of a center vertex (the \emph{hub}) and a
ring of $n$ peripheral vertices, each connected to the hub; see
Figure~\ref{f-wheel}.  It has total $n+1$ nodes.  Its automorphism
group is isomorphic to the dihedral group $D_n$ with order $2n$.  The
transition probability matrix has the structure
\begin{equation}
P=\left[\begin{array}{cccccc}
1-n p & p & p & \ldots & p & p \\
p & 1-p-2q & q & \ldots & 0 & q \\
p & q & 1-p-2q & \ldots & 0 & 0 \\
\vdots & \vdots & \vdots & \ddots & \vdots & \vdots \\
p & 0 & 0 & \ldots & 1-p-2q & q \\
p & q & 0 & \ldots & q & 1 -p-2q 
\end{array}\right],
\end{equation} 
where $p$ and $q$ are the transition probabilities between the hub and 
each peripheral vertex, and between adjacent peripheral vertices, 
respectively.

For this structure, the block-diagonalizing transformation is given by
\[
T = \Diag(1,\mathcal{F}_n), \qquad 
[\mathcal{F}_n]_{jk} = \frac{1}{\sqrt{n}}
e^\frac{2 \pi \imath (j-1)(k-1)}{n}
\]
where $\mathcal{F}_n$ is the unitary Fourier matrix of size 
$n\times n$. 
As a consequence, the matrix $T^{-1} P T$ is block diagonal with 
a $2\times 2$ matrix and $n-1$ scalars on its diagonal, given by
\[
\left[\begin{array}{cc}
1-np & \sqrt{n} p \\
\sqrt{n} p & 1 - p
\end{array}\right] \\
\]
and
\[
1 - p  + (\omega_n^k + \omega_n^{-k} - 2 ) \cdot q, 
\qquad k=1,\ldots,n-1
\]
where $\omega_n= e^\frac{2 \pi \imath}{n}$ is an elementary $n$-th
root of unity.
The $2\times 2$ block is $B_1$, which contains eigenvalues of the
orbit chain under $D_n$ (it has only two orbits). 

With the above decomposition, we obtain the optimal solution to the
FMMC problem in closed form
\[
p^\star = \frac{1}{n}, \qquad 
q_\star =  \frac{1-\frac{1}{n}}{2-\cos\frac{2\pi}{n}
-\cos\frac{2\lfloor n/2\rfloor\pi}{n}}.
\]
The optimal value of the SLEM is
\[
\mu^{\star} = \left(1-\frac{1}{n}\right)
  \frac{\cos\frac{2\pi}{n}-\cos\frac{2\lfloor n/2\rfloor\pi}{n}}
  {2-\cos\frac{2\pi}{n}-\cos\frac{2\lfloor n/2\rfloor\pi}{n}}.
\]
Compared with the optimal solution for the cycle graph
in~(\ref{e-cycle-prob}) and~(\ref{e-cycle-slem}), we see an extra
factor of $1-1/n$ in both the SLEM and the transition probability
between peripheral vertices.
This is exactly the factor improved by adding the central hub 
over the pure $n$-cycle case.

The wheel graph is an example for which the block diagonalization
technique works out nicely, while the orbit theory leads to much less 
reduction. 
Although there are only two orbits under the full automorphism group,
any orbit graph that has a fixed peripheral vertex
will have at least $(n+1)/2$ orbits (the corresponding symmetry is
the reflection through that vertex).  

\subsubsection{$K_n$-$K_n$}
\label{s-Kn-Kn-2}
We did careful symmetry analysis for the graph $K_n$-$K_n$
in~\S\ref{s-Kn-Kn}; see Figure~\ref{f-KnKn}.
The transition probability matrix on this graph has the structure 
\[
P = \left[ \begin{array}{cccc}
C & p_1 \ones & 0 & 0 \\
p_1 \ones^T & 1-p_0-(n-1)p_1 & p_0 & 0 \\
0 & p_0 & 1-p_0-(n-1)p_1 & p_1 \ones^T \\
0 & 0 & p_1 \ones & C  \end{array} \right]
\]
where $C$ is a circulant matrix
\[
C = (1-p_1-(n-3)p_2)I_{n-1} + p_2 \ones_{(n-1)\times(n-1)}.
\]

Since circulant matrices are diagonalized by Fourier matrices, we
first use the transformation matrix
\[
T_1 = \left[ \begin{array}{cccc}
\mathcal F_{n-1} & 0 & 0 & 0\\
0 & 1 & 0 & 0 \\
0 & 0 & 1 & 0 \\
0 & 0 & 0 & \mathcal F_{n-1} \end{array} \right]
\] 
where $\mathcal F_{n-1}$ is the unitary Fourier matrix of
dimension $n-1$. 
This corresponds to block diagonalization using the symmetry group
$S_{n-1}\times S_{n-1}$, which is a subgroup of
$\mbox{Aut}(\mbox{$K_n$-$K_n$})$.  
The matrix $T_1^{-1} P T_1$ has diagonal blocks
\[
B_1'=\left[ \begin{array}{cccc}
1-p_1 & \sqrt{n-1}p_1 & 0 & 0\\
\sqrt{n-1}p_1 & 1-p_0-(n-1)p_1 & p_0 & 0\\
0 & p_0 & 1-p_0-(n-1)p_1 & \sqrt{n-1}p_1\\
0 & 0 & \sqrt{n-1}p_1 & 1-p_1 \end{array} \right] 
\]
and
\begin{equation}\label{e-KnKn-2n-4}
I_{2n-4} \otimes (1-p_1-(n-1)p_2).
\end{equation}
From this we know that $P$ has an eigenvalue $1-p_1-(n-1)p_2$ with
multiplicity $2n-4$, and the remaining four eigenvalues are the
eigenvalues of the above $4\times 4$ block $B_1'$.  The block $B_1'$
corresponds to the orbit chain under the symmetry group
$H=S_{n-1}\times S_{n-1}$.  More precisely, $B_1'=\Pi^{1/2} P_H
\Pi^{-1/2}$, where $\Pi=\Diag(\pi_H)$, $P_H$ and $\pi_H$ are the
transition probability matrix and stationary distribution of the orbit
chain shown in Figure~\ref{f-KnKn:c}, respectively.

Exploring the full automorphism group of $K_n$-$K_n$, we can further
block diagonalize $B_1'$. 
Let
\[
T = T_1 \left[ \begin{array}{ccc}
I_{n-2} & & \\ & T_2 & \\ & & I_{n-2} \end{array} \right], \qquad
T_2 = \frac{1}{\sqrt{2}} \left[ \begin{array}{ccrr}
1 & 0 & 0 & 1\\ 0 & 1 & 1 & 0 \\ 0 & 1 & -1 & 0\\ 1 & 0 & 0 & -1
\end{array} \right].
\]
The $4\times 4$ block $B_1'$ is decomposed into
\[
\left[ \begin{array}{cc}
1-p_1 & \sqrt{n-1}p_1 \\
\sqrt{n-1}p_1 & 1-(n-1)p_1 
\end{array} \right], 
\qquad
\left[ \begin{array}{cc}
1-2p_0-(n-1)p_1 & \sqrt{n-1}p_1\\
\sqrt{n-1}p_1 & 1-p_1 \end{array} \right] 
\]
The first block is $B_1$, which has eigenvalues~$1$ and $1-np_1$.
By Theorem~\ref{t-orbit-diag-connection}, $B_1$ is related to the
orbit chain under $\mbox{Aut}(\mbox{$K_n$-$K_n$})$ 
(see Figure~\ref{f-KnKn:b}) by the
equation~(\ref{e-orbit-diag-connection}). 
The second $2\times 2$ block has eigenvalues
\[
1-p_0-(1/2)np_1\pm\sqrt{(p_0+(1/2)np_1)^2-2p_0p_1}.
\]
These are the eigenvalues contained in the orbit chain of
Figure~\ref{f-KnKn:c} but not in Figure~\ref{f-KnKn:b}.

In summary, the distinct eigenvalues of the Markov chain on
$K_n$-$K_n$ are
\[
1,\quad 1-np_1, \quad
1-p_0-(1/2)np_1\pm\sqrt{(p_0+(1/2)np_1)^2-2p_0p_1},
\quad 1-p_1-(n-1)p_2
\]
where the last one has multiplicity $2n-4$, and all the rest have
multiplicity~1.  To solve the FMMC problem, we still need to solve the
SDP~(\ref{e-fmmc-block-diag}).  There are three blocks of matrix
inequality constraints, with sizes $2$, $2$, $1$, respectively.  Note
that the total size is $5$, which is exactly the size of the single
matrix inequality in the SDP~(\ref{e-fmmc-rev}) when we used the orbit
theory to do symmetry reduction. As we mentioned before, the huge
reduction for $K_n$-$K_n$ is due to the fact that it has an
irreducible representation with high dimension $2n-4$ and
multiplicity~$1$ (see \cite[Proposition 2.4]{BDPX:05}).  In the
decomposition~(\ref{e-block-diag}), this means a block of size~$1$
repeated $2n-4$ times; see equation~(\ref{e-KnKn-2n-4}).

Since now the problem has been reduced to something much more
tractable, we can even obtain an analytic expression for the optimal
transition probabilities. The optimal solution for the $K_n$-$K_n$
graph (for $n \geq 2$) is given by:
\[
p_0^\star = (\sqrt{2}-1) \frac{n+\sqrt{2}-2}{n+2-2\sqrt{2}}, \qquad
p_1^\star = \frac{2 - \sqrt{2}}{n+2-2\sqrt{2}}, \qquad
p_2^\star = \frac{n-\sqrt{2}}{(n-1)(n+2-2\sqrt{2})}.
\]
The corresponding optimal convergence rate is 
\[
\mu^\star = \frac{n-4+2\sqrt{2}}{n+2-2\sqrt{2}}.
\]
For large $n$, we have 
$\mu^\star=1-\frac{6 - 4 \sqrt{2}}{n} + O\left(\frac{1}{n^2}\right)$.
This is quite close to the SLEM of a suboptimal construction with 
transition probabilities 
\begin{equation}\label{e-sub-optimal}
p_0 = \frac{1}{2}, \qquad p_1=p_2 = \frac{1}{2(n-1)}.
\end{equation}
As shown in \cite{BDPX:05}, the corresponding SLEM is of the order 
$\mu = 1 - \frac{1}{3n} + O\left(\frac{1}{n^2}\right)$; 
here we have $6 - 4 \sqrt{2} \approx 0.3431$. 
The limiting value of the optimal
transition probability between the two clusters is $\sqrt{2}-1 \approx
0.4142$.

\subsubsection{Complete binary trees}

Since the automorphism groups of the complete binary trees $\mathcal
T_n$ are given recursively (see~\S\ref{s-tree-orbit}), it is also
convenient to write the transition probability matrices in a recursive
form. 
We start from the bottom by considering the last level of branches. 
If we cut-off the rest of the tree, the last level has three nodes and
two edges with the transition probability matrix
\begin{equation}
P_n = 
\left[\begin{array}{ccc}
1-2 p_n & p_n & p_n \\
p_n & 1-p_n & 0 \\
p_n & 0 & 1-p_n
\end{array}\right].
\end{equation}
For the tree with $n$ levels $\mathcal T_n$, the transition matrix $P_1$
can be computed from the recursion
\begin{equation}
P_{k-1} = 
\left[\begin{array}{ccc}
1-2 p_{k-1} & p_{k-1}e_k^T & p_{k-1}e_k^T\\
p_{k-1}e_k & P_k - p_{k-1} e_k e_k^T & 0\\ 
p_{k-1}e_k & 0 & P_k - p_{k-1} e_k e_k^T 
\end{array}\right], \qquad k=n,n-1\ldots,2
\end{equation}
where $e_k=[1~0\ldots~0]$, a unit vector in $\reals^{t_k}$ with 
$t_k = 2^{k+1}-1$. 

The coordinate transformations are also best written in recursive
form. 
Let 
\[
T_n \Diag(1,\mathcal{F}_2), \qquad 
\mathcal{F}_2 = 
\frac{1}{\sqrt{2}}
\left[\begin{array}{cc}
1 & 1 \\
1 & -1
\end{array}
\right],
\] 
and define the matrices
\[
T_{k-1} = \Diag(1,\mathcal{F}_2 \otimes T_k),
\qquad k=n,n-1,\ldots,2.
\]
It is clear that all the $T_k$ are orthogonal.
It is easy to verify that $T_n$ block-diagonalizes $P_n$ 
\[
T_n^T P_n T_n = \left[ \begin{array}{ccc}
1-2p_n & \sqrt{2}p_n & 0 \\
\sqrt{2}p_n & 1-p_n & 0\\
0 & 0 & 1-p_n \end{array} \right] .
\]
In fact $T_k$ block-diagonalizes $P_k$, and the transformed matrices 
can be obtained recursively
\[
T_{k-1}^T P_{k-1} T_{k-1} = \left[ \begin{array}{ccc}
1-2p_{k-1} & \sqrt{2}p_{k-1} e_k^T & 0\\ \sqrt{2}p_{k-1} e_k & 
T_k^T P_k T_k - p_{k-1} e_k e_k^T & 0\\ 
0 & 0 &  T_k^T P_k T_k - p_{k-1} e_k e_k^T 
\end{array} \right]
\]
for $k=n,n-1,\ldots,2$.

The matrix $T_1^T P_1 T_1$ has a very special structure. 
It has $n+1$ distinct blocks, each with size $1,\ldots,n+1$,
respectively.  
Order these blocks with increasing sizes as $B_1,B_2,\ldots,B_{n+1}$. 
The largest block of size $n+1$ is
\[
B_{n+1} = \left[ \begin{array}{cccccc}
1\!-\!2p_1 & \sqrt{2}p_1 & & & & \\
\sqrt{2}p_1 & 1\!-\!p_1\!-\!2p_2 & \sqrt{2}p_2 & & & \\
 & \sqrt{2}p_2 & 1\!-\!p_2\!-\!2p_3 & \sqrt{2} p_3 & & \\[2ex]
 & & \ddots & \ddots & \ddots & \\[2ex]
 & & & \sqrt{2}p_{n\!-\!1} & 1\!-\!p_{n\!-\!1}\!-\!2p_n & \sqrt{2}p_n \\
 & & & & \sqrt{2}p_n & 1\!-\!p_n
\end{array} \right].
\]
The matrix $B_n$ is the submatrix of~$B_{n+1}$ by removing its first
row and column. 
The matrix $B_{n-1}$ is the submatrix of~$B_{n+1}$ by removing its
first two rows and first two columns, and so on.
The matrix $B_1$ is just the scalar $1-p_n$.  
The matrix $B_{n+1}$ only appears once and it is related 
by~(\ref{e-orbit-diag-connection}) to the orbit chain in 
Figure~\ref{f-binary-tree:a} (for this example we use
$B_{n+1}$ instead of $B_1$ for notational convenience).
The eigenvalues of $B_{n+1}$ appear in $\mathcal T_n$ with
multiplicity one.   
For $k=1,\ldots,n$, the block $B_k$ is repeated $2^{n-k}$
times.  
These blocks, in a recursive form, contain additional eigenvalues of
$\mathcal T_n$, and the numbers of their occurrences reveal the
multiplicities of the eigenvalues. 

More specifically, we note that the orbit chain under the full
automorphism group has only one fixed point --- the root vertex 
(see Figure~\ref{f-binary-tree:a}).
We consider next the orbit chain that has a fixed point in the first
level of child vertices (the other child vertex in the same level 
is also fixed).
This is the orbit graph in Figure~\ref{f-binary-tree:b}, which has
$2n+1$ vertices. 
The matrix $B_n$ contains exactly the~$n$
eigenvalues that appear in this orbit chain but not in the one of
Figure~\ref{f-binary-tree:a}. 
These~$n$ eigenvalues each has multiplicity $2^{n-n}=1$ 
in $\mathcal T_n$.  
Then we consider the orbit chain that has a fixed point in the second
level of child vertices (it also must have a fixed point in the
previous level).
This is the orbit graph in Figure~\ref{f-binary-tree:c}, which has
$3n$ vertices.
The matrix $B_{n-1}$ contains exactly the $n-1$
eigenvalues that appear in this orbit chain but not in the previous
one. 
These $n-1$ eigenvalues each has multiplicity $2^{n-(n-1)}=2$. 
In general, for $k=1,\ldots,n$, 
the size of the orbit chain that has a fixed point in the
$k$-th level of child vertices is
\[
(n+1) + n + \cdots + (n+1-k) 
\]
(it must have a fixed point in all previous levels).
Compared with the orbit chain of $(k-1)$-th level, 
the orbit chain of $k$-th level contains additional $n+1-k$
eigenvalues.
These are precisely the eigenvalues of the matrix $B_{n+1-k}$, 
and they all appear in $\mathcal T_n$ with multiplicity 
$2^{n-(n+1-k)}=2^{k-1}$.

Because of the special structure of $B_1,\ldots,B_{n+1}$, we have the
following eigenvalue interlacing result 
(\eg, \cite[Theorem 4.3.8]{HoJ:85})
\[
\lambda_{k+1}(B_{k+1}) \leq \lambda_k(B_k) \leq \lambda_k(B_{k+1})
\leq \lambda_{k-1}(B_k) \leq \cdots \leq \lambda_2(B_k) \leq 
\lambda_2(B_{k+1}) \leq \lambda_1(B_k) \leq \lambda_1(B_{k+1}) 
\]
for $k=1,\ldots,n$.
Thus for the FMMC problem, we only need to consider the
two blocks $B_{n+1}$ and $B_n$ (note that $\lambda_1(B_{n+1})=1$).  
In other words, we only need to consider the orbit chain with $2n+1$
vertices in Figure~\ref{f-binary-tree:b}. 
This is a further simplification over the method based on orbit theory.  

We conjecture that the optimal transition probabilities are
\[
p^\star_k =  \frac{1}{3} \left( 1 - \left(- \frac{1}{2} \right)^k \right), 
\qquad k = 1,\ldots,n.
\]
Notice that these probabilities \emph{do not} depend explicitly
on $n$, and so they coincide for any two binary trees, regardless of
the height.  With increasing $k$, the limiting optimal values
oscillate around and converge to $1/3$.

\subsubsection{An example of Ron Graham}
\label{s-graham-example}

\begin{figure}
\begin{center}
\includegraphics[width=0.8\textwidth]{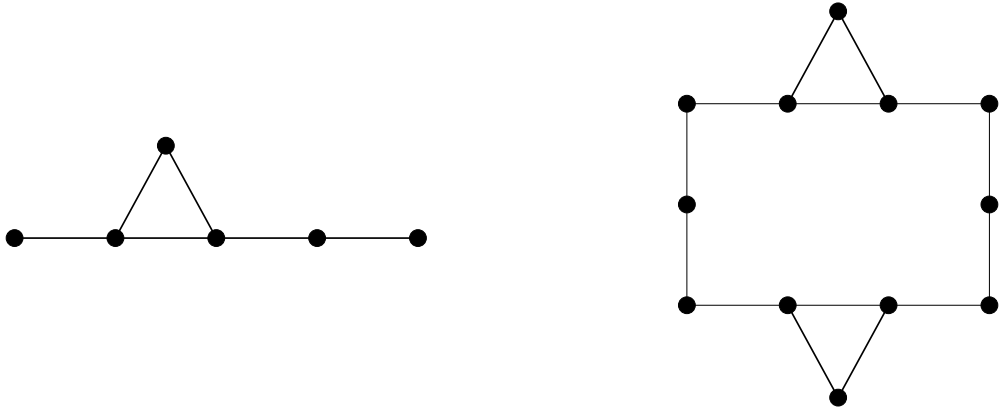}
\caption{Left: the simplest graph with no symmetry.
Right: two copies joined head-to-tail.}
\label{f-graham-graph}
\end{center}
\end{figure}

\begin{figure}
\begin{center}
\psfrag{p1}[bc]{$p_1$}
\psfrag{p2}[bc]{$p_2$}
\psfrag{p3}[br]{$p_3$}
\psfrag{p3r}[bl]{$p_3$}
\psfrag{p4}[bl]{$p_4$}
\psfrag{2p3}[tr]{$2p_3$}
\psfrag{2p4}[tc]{$2p_4$}
\includegraphics[width=0.8\textwidth]{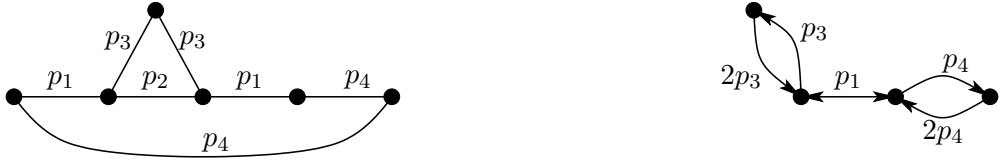}
\caption{Left: orbit graph with $C_n$ symmetry.
Right: orbit graph with $D_n$ symmetry.}
\label{f-graham-orbits}
\end{center}
\end{figure}

We finish this section with an example raised by Ron Graham.
Consider the simplest graph with no symmetry
(Figure~\ref{f-graham-graph}, left). 
Take~$n$ copies of this six vertex graph and join them, head to tail,
in a cycle.
By construction, this $6n$ vertex graph certainly has $C_n$ symmetry.
Careful examination reveals that the automorphism group is isomorphic to the dihedral group $D_n$ (with order $2n$). 
The construction actually brings symmetry under reflections in addition to rotations (Figure~\ref{f-graham-graph}, right). 
The orbit graphs under $C_n$ and $D_n$ are shown in Figure~\ref{f-graham-orbits}.

Although the automorphism group of this graph (with $6n$ vertices) is isomorphic to the ones of $n$-cycles (Figure~\ref{f-cycle}) and wheels (Figure~\ref{f-wheel}), finding the symmetry-adapted basis for block-diagonalization is a bit more involved. 
This is due to the different types of orbits we have for this graph.
The details of block-diagonalizing this type of graphs is described in \cite[\S3.1]{FaS:92}. 
The diagonal blocks of the resulting matrix all have sizes no larger than $6\times 6$. 
Numerical experiments show that for $n\geq 3$, the fastest mixing chain seems to satisfy
\[
p^\star_1=p^\star_4=\frac{1}{2}, \qquad p^\star_2+p^\star_3=\frac{1}{2}.
\] 
Intuitively, this $6n$ vertex graph is the same as modifying a $5n$ vertex cycle by adding a triangular bump (with an additional vertex) for every $5$ vertices.
Recall that for a pure cycle, we have to use a transition probability that is slightly less than $1/2$ to achieve fastest mixing; see equation~(\ref{e-cycle-prob}). 
Here because of the added bumps, it seems optimal to assign transition probability $1/2$ to every edge on the cycle ($p^\star_1$ and $p^\star_4$), except for edges being part of a bump.
For the bumps, the probability $1/2$ is shared between the original edge on the cycle ($p^\star_2$) and the edge connecting to the bump points ($p^\star_3$).   
Moreover, we observe that as $n$ increases, $p^\star_3$ gets smaller and $p^\star_2$ gets closer to $1/2$. 
So for large $n$, the added bump vertices seem to be ignored, with very small probability to be reached; but once it is reached, it will staying there with high probability.

\section{Conclusions}
\label{s-conclusions}

We have shown that exploiting graph symmetry can lead to significant
reduction in both the number of variables and the size of matrices, in
solving the FMMC problem.  For special classes of graphs such as
edge-transitive and distance-transitive graphs, symmetry reduction
leads to closed form solutions in terms of the eigenvalues of the
Laplacian matrix or the intersection matrix.  For more general graphs,
we gave two symmetry reduction methods, based on orbit theory and
block diagonalization, respectively.

The method based on orbit theory is very intuitive, but the
construction of ``good'' orbit chains can be of more art than
technique.  The method of block diagonalization can be mostly
automated once the irreducible representations of the automorphism
groups are generated (for small graphs, they can be generated using
software for computational discrete algebra such as GAP \cite{GAP}).
These two approaches have an interesting connection: orbit theory
gives nice interpretation of the diagonal blocks, while the block
diagonalization approach offers theoretical insights about the
construction of the orbit chains.

The symmetry reduction method developed in this paper can be very
useful in many combinatorial optimization problems where the graph has
rich symmetry properties, in particular, problems that can be
formulated as or approximated by SDP or eigenvalue optimization
problems involving weighted Laplacian matrices (\eg,
\cite{MoP:93,Goe:97}).  In addition to the reduction of problem size,
other advantages of symmetry exploitation includes degeneracy removal,
better conditioning and reliability \cite{GatermannParrilo}.

There is still much to do in understanding how to exploit symmetry in
semidefinite programming.  The techniques presented in this paper
requires a good understanding of orbit theory, group representation
theory and interior-point methods for SDP.  It is of practical
importance to develop general purpose methods that can automatically
detect symmetries (\eg, the code \texttt{nauty} \cite{nauty} for
graph automorphisms), and then exploit them in computations.  A good
model here is general purpose (but heuristic) methods for exploiting
sparsity in numerical linear algebra, where symbolic operations on
graphs (\eg, minimum degree permutation) reduce fill-ins in numerical
factorization (\eg, \cite{GeL:81}).  As a result of this work, even
very large sparse optimization problems are now routinely solved by
users who are not experts in sparse matrix methods.  For exploiting
symmetry in SDP, the challenges include the development of fast
methods to detect large symmetry groups (for computational purposes,
it often suffices to recognize parts of the symmetries), and the
integration of algebraic methods (\eg, orbit theory and group
representations) and numerical algorithms (\eg, interior-point
methods).

\bibliographystyle{alpha}
\bibliography{fmmc_symmetry}

\begin{thebibliography}{KOMK01}

\bibitem[AV05]{AbV:05}
M.~Ab\'ert and B.~Vir\'ag.
\newblock Dimension and randomness in groups acting on rooted trees.
\newblock {\em Journal of the American Mathematical Society}, 18(1):157--192,
  2005.

\bibitem[BCN89]{BCN:89}
A.~E. Brouwer, A.~M. Cohen, and A.~Neumaier.
\newblock {\em Distance-Regular Graphs}.
\newblock Springer-Verlag, Berlin, 1989.

\bibitem[BDPX05]{BDPX:05}
S.~Boyd, P.~Diaconis, P.~A. Parrilo, and L.~Xiao.
\newblock Symmetry analysis of reversible {M}arkov chains.
\newblock {\em Internet Mathematics}, 2(1):31--71, 2005.

\bibitem[BDSX06]{BDSX:06}
S.~Boyd, P.~Diaconis, J.~Sun, and L.~Xiao.
\newblock Fastest mixing {M}arkov chain on a path.
\newblock {\em The American Mathematical Monthly}, 113(1):70--74, January 2006.

\bibitem[BDX04]{BDX:04}
S.~Boyd, P.~Diaconis, and L.~Xiao.
\newblock Fastest mixing {Markov} chain on a graph.
\newblock {\em SIAM Review}, 46(4):667--689, 2004.

\bibitem[Big74]{Biggs}
N.~Biggs.
\newblock {\em Algebraic Graph Theory}.
\newblock Cambridge University Press, 1974.

\bibitem[BM03]{BuM:03}
S.~Burer and R.~D.~C. Monteiro.
\newblock A nonlinear programming algorithm for solving semidefinite programs
  via low-rank factorization.
\newblock {\em Mathematical Programming, Series B}, 95:329--357, 2003.

\bibitem[Br{\'e}99]{Bre:99}
P.~Br{\'e}maud.
\newblock {\em Markov Chains, Gibbs Fields, Monte Carlo Simulation and Queues}.
\newblock Texts in Applied Mathematics. Springer-Verlag, Berlin-Heidelberg,
  1999.

\bibitem[BTN01]{BTN:01}
A.~Ben-Tal and A.~Nemirovski.
\newblock {\em Lectures on Modern Convex Optimization, Analysis, Algorithms,
  and Engineering Applications}.
\newblock MPS/SIAM Series on Optimization. SIAM, 2001.

\bibitem[BV04]{BoV:04}
S.~Boyd and L.~Vandenberghe.
\newblock {\em Convex Optimization}.
\newblock Cambridge University Press, 2004.
\newblock Available at \verb+http://www.stanford.edu/~boyd/cvxbook.html+.

\bibitem[BYZ00]{BYZ:00}
S.~Benson, Y.~Ye, and X.~Zhang.
\newblock Solving large-scale sparse semidefinite programs for combinatorial
  optimization.
\newblock {\em SIAM Journal Optimization}, 10:443--461, 2000.

\bibitem[Chu97]{FanChung}
F.~R.~K. Chung.
\newblock {\em Spectral Graph Theory}.
\newblock Number~92 in CBMS Regional Conference Series in Mathematics. American
  Mathematical Society, 1997.

\bibitem[CLP03]{CLP:03}
R.~Cogill, S.~Lall, and P.~A. Parrilo.
\newblock On structured semidefinite programs for the control of symmetric
  systems.
\newblock In {\em Proceedings of 41st Allerton Conference on Communication,
  Control, and Computing}, pages 1536--1545, Monticello, IL, October 2003.

\bibitem[DCS80]{DCS:80}
M.~Doob, D.~Cvetkovic, and H.~Sachs.
\newblock {\em Spectra of Graphs: Theory and Application}.
\newblock Academic Press, New York, 1980.

\bibitem[Dia88]{Dia:88}
P.~Diaconis.
\newblock {\em Group Representations in Probability and Statistics}.
\newblock IMS, Hayward, CA, 1988.

\bibitem[dKPS07]{dPS:07}
E.~de~Klerk, D.~V. Pasechnik, and A.~Schrijver.
\newblock Reduction of symmetric semidefinite programs using the regular
  $*$-representation.
\newblock {\em Mathematical Programming, Series B}, 109:613--624, 2007.

\bibitem[DR07]{DuR:07}
I.~Dukanovic and F.~Rendl.
\newblock Semidefinite programming relaxations for graph coloring and maximal
  clique problems.
\newblock {\em Mathematical Programming, Series B}, 109:345--365, 2007.

\bibitem[DS91]{DiS:91}
P.~Diaconis and D.~Stroock.
\newblock Geometric bounds for eigenvalues of {M}arkov chains.
\newblock {\em The Annals of Applied Probability}, 1(1):36--61, 1991.

\bibitem[DSC93]{DSC:93}
P.~Diaconis and L.~Saloff-Coste.
\newblock Comparison theorems for reversible {M}arkov chains.
\newblock {\em Ann. Appl. Probab.}, 3:696--730, 1993.

\bibitem[DSC06]{DSC:06}
P.~Diaconis and L.~Saloff-Coste.
\newblock Separation cut-offs for birth and death chains.
\newblock Submitted to \emph{Annals of Applied Probabilities}, 2006.

\bibitem[ER63]{ErR:63}
P.~Erd{\H o}s and R\'enyi.
\newblock Asymmetric graphs.
\newblock {\em Acta Math. Acad. Sci. Hungar.}, 14:295--315, 1963.

\bibitem[FS92]{FaS:92}
A.~F\"assler and E.~Stiefel.
\newblock {\em Group Theoretical Methods and Their Applications}.
\newblock Birkh\"auser, Boston, 1992.

\bibitem[Gat00]{Gat:00}
K.~Gatermann.
\newblock {\em Computer Algebra Methods for Equivariant Dynamical Systems},
  volume 1728 of {\em Lecture Notes in Mathematics}.
\newblock Springer-Verlag, 2000.

\bibitem[GL81]{GeL:81}
A.~George and J.~Liu.
\newblock {\em Computer Solution of Large Sparse Positive Definite Systems}.
\newblock Prentice Hall, Englewood Cliffs, NJ, 1981.

\bibitem[GL96]{GoV:96}
G.~H. Golub and C.~F.~Van Loan.
\newblock {\em Matrix Computations}.
\newblock The John Hopkins University Press, Baltimore, 3rd edition, 1996.

\bibitem[GO06]{GaO:06}
K.~K. Gade and M.~L. Overton.
\newblock Optimizing the asymptotic convergence rate of the
  {Diaconis-Holmes-Neal} sampler.
\newblock To appear in \emph{Advances in Applied Mathematics}, 2006.

\bibitem[Goe97]{Goe:97}
M.~X. Goemans.
\newblock Semidefinite programming in combinatorial optimization.
\newblock {\em Mathematical Programming}, 79:143--161, 1997.

\bibitem[GP04]{GatermannParrilo}
K.~Gatermann and P.~A. Parrilo.
\newblock Symmetry groups, semidefinite programs, and sums of squares.
\newblock {\em Journal of Pure and Appl. Algebra}, 192:95--128, 2004.

\bibitem[GR01]{GoR:01}
C.~Godsil and G.~Royle.
\newblock {\em Algebraic Graph Theory}, volume 207 of {\em Graduate Texts in
  Mathematics}.
\newblock Springer, 2001.

\bibitem[Gra81]{Gra:81}
A.~Graham.
\newblock {\em Kronecker Products and Matrix Calculus with Applications}.
\newblock Ellis Horwoods Ltd., Chichester, UK, 1981.

\bibitem[gro05]{GAP}
The~GAP group.
\newblock {GAP} - groups, algorithms, programming - a system for computational
  discrete algebra, version 4.4.6, 2005.
\newblock \texttt{http://www.gap-system.org}.

\bibitem[GSS88]{GSS:88}
M.~Golubitsky, I.~Stewart, and D.~G. Schaeffer.
\newblock {\em Singularities and Groups in Bifurcation Theory II}, volume~69 of
  {\em Applied Mathematical Sciences}.
\newblock Springer, New York, 1988.

\bibitem[HJ85]{HoJ:85}
R.~A. Horn and C.~A. Johnson.
\newblock {\em Matrix Analysis}.
\newblock Cambridge University Press, 1985.

\bibitem[HOY03]{HOY:03}
B.~Han, M.~L. Overton, and T.~P.-Y. Yu.
\newblock Design of {Hermite} subdivision schemes aided by spectral radius
  optimization.
\newblock {\em SIAM Journal on Matrix Analysis and Applications}, 25:80--104,
  2003.

\bibitem[HR00]{HeR:00}
C.~Helmberg and F.~Rendl.
\newblock A spectral bundle method for semidefinite programming.
\newblock {\em SIAM Journal on Optimization}, 10(3):673--696, 2000.

\bibitem[JK81]{JaK:81}
G.~D. James and A.~Kerber.
\newblock {\em The Representation Theory of the Symmetric Group}.
\newblock Addison-Wesley, Reading, Massachusetts, 1981.

\bibitem[KOMK01]{KOMK:01}
Y.~Kanno, M.~Ohsaki, K.~Murota, and N.~Katoh.
\newblock Group symmetry in interior-point methods for semidefinite
  programming.
\newblock {\em Optimization and Engineering}, 2:293--320, 2001.

\bibitem[Lau07]{Lau:07}
M.~Laurent.
\newblock Strengthend semidefinite programming bounds for codes.
\newblock {\em Mathematical Programming, Series B}, 109:239--261, 2007.

\bibitem[LNM04]{LNM:04}
Z.~Lu, A.~Nemirovski, and R.~D.~C. Monteiro.
\newblock Large-scale semidefinite programming via saddle point mirror-prox
  algorithm.
\newblock {\em Submitted to \emph{Mathematical Programming}}, 2004.

\bibitem[Mar03]{Mar:03}
F.~Margot.
\newblock Exploiting orbits in symmetric {ILP}.
\newblock {\em Mathematical Programming, Series B}, 98:3--21, 2003.

\bibitem[McK03]{nauty}
B.D. McKay.
\newblock {\em \texttt{nauty} User's guide (Version 2.2)}.
\newblock Australian National University, 2003.
\newblock Available from \texttt{http://cs.anu.edu.au/\~{}bdm/nauty/}.

\bibitem[Mer94]{Mer:94}
R.~Merris.
\newblock {L}aplacian matrices of graphs: a survey.
\newblock {\em Linear Algebra and Its Applications}, 197:143--176, 1994.

\bibitem[Moh97]{Moh:97}
B.~Mohar.
\newblock Some applications of {L}aplace eigenvalues of graphs.
\newblock In G.~Hahn and G.~Sabidussi, editors, {\em Graph Symmetry: Algebraic
  Methods and Applications}, NATO ASI Ser. C 497, pages 225--275. Kluwer, 1997.

\bibitem[MP93]{MoP:93}
B.~Mohar and S.~Poljak.
\newblock Eigenvalues in combinatorial optimization.
\newblock In R.~A. Brualdi, S.~Friedland, and V.~Klee, editors, {\em
  Combinatorial and Graph-Theoretical Problems in Linear Algebra}, volume~50 of
  {\em IMA Volumes in Mathematics and Its Applications}, pages 107--151.
  Springer-Verlag, 1993.

\bibitem[MR99]{MaR:99}
J.~E. Marsden and T.~Ratiu.
\newblock {\em Introduction to Mechanics and Symmetry}, volume~17 of {\em Texts
  in Applied Mathematics}.
\newblock Springer-Verlag, 2nd edition, 1999.

\bibitem[Nem04]{Nem:04}
A.~Nemirovski.
\newblock Prox-method with rate of convergence {$O(1/t)$} for variational
  inequalities with {L}ipschitz continuous monotone operators and smooth
  convex-concave saddle point problems.
\newblock {\em SIAM Journal on Optimization}, 15(1):229--251, 2004.

\bibitem[Nes05]{Nes:05}
Y.~Nesterov.
\newblock Smooth minimization of non-smooth functions.
\newblock {\em Mathematical Programming}, 103:127--152, 2005.

\bibitem[NN94]{NeN:94}
Y.~Nesterov and A.~Nemirovskii.
\newblock {\em Interior-Point Polynomial Algorithms in Convex Programming}.
\newblock SIAM Studies in Applied Mathematics. SIAM, 1994.

\bibitem[OOR04]{OOR:04}
R.~C. Orellana, M.~E. Orrison, and D.~N. Rockmore.
\newblock Rooted trees and iterated weath products of cyclic groups.
\newblock {\em Advances in Applied Mathematics}, 33(3):531--547, 2004.

\bibitem[Ove92]{Ove:92}
M.~L. Overton.
\newblock Large-scale optimization of eigenvalues.
\newblock {\em SIAM Journal on Optimization}, 2:88--120, 1992.

\bibitem[Par00]{PhD:Parrilo}
P.~A. Parrilo.
\newblock {\em Structured semidefinite programs and semialgebraic geometry
  methods in robustness and optimization}.
\newblock PhD thesis, California Institute of Technology, May 2000.
\newblock Available at
  \texttt{http://resolver.caltech.edu/CaltechETD:etd-05062004-055516}.

\bibitem[Par03]{Par:03}
P.~A. Parrilo.
\newblock Semidefinite programming relaxations for semialgebraic problems.
\newblock {\em Mathematical Programming}, 96:293 -- 320, 2003.

\bibitem[PS03]{ParriloSturmfels}
P.~A. Parrilo and B.~Sturmfels.
\newblock Minimizing polynomial functions.
\newblock In S.~Basu and L.~Gonzalez-Vega, editors, {\em Algorithmic and
  quantitative real algebraic geometry}, volume~60 of {\em DIMACS Series in
  Discrete Mathematics and Theoretical Computer Science}, pages 83--99. AMS,
  2003.

\bibitem[Roc05]{Roc:05}
S.~Roch.
\newblock Bounding fastest mixing.
\newblock {\em Electronic Communications in Probability}, 10:282--296, 2005.

\bibitem[Saa92]{Saa:92}
Y.~Saad.
\newblock {\em Numerical Methods for Large Eigenvalue Problems}.
\newblock Manchester University Press, Manchester, UK, 1992.

\bibitem[Sal06]{Sal:06}
J.~Saltzman.
\newblock {\em A generalization of spectral analysis for discrete data using
  {Markov} chains}.
\newblock PhD thesis, Department of Statistics, Stanford University, 2006.

\bibitem[Ser77]{Ser:77}
J.-P. Serre.
\newblock {\em Linear Representations of Finite Groups}.
\newblock Springer-Verlag, New York, 1977.

\bibitem[Stu99]{Stu:99}
J.~F. Sturm.
\newblock Using {S}e{D}u{M}i 1.02, a {MATLAB} toolbox for optimization over
  symmetric cones.
\newblock {\em Optimization Methods and Software}, 11-12:625--653, 1999.
\newblock Special issue on Interior Point Methods (CD supplement with
  software).

\bibitem[TTT99]{SDPT3}
K.~C. Toh, M.~J. Todd, and R.~H. Tutuncu.
\newblock {SDPT3} --- a {Matlab} software package for semidefinite programming.
\newblock {\em Optimization Methods and Software}, 11:545--581, 1999.

\bibitem[VB96]{VaB:96}
L.~Vandenberghe and S.~Boyd.
\newblock Semidefinite programming.
\newblock {\em SIAM Review}, 38(1):49--95, 1996.

\bibitem[Wor94]{Wor:94}
P.~Worfolk.
\newblock Zeros of equivariant vector fields: Algorithms for an invariant
  approach.
\newblock {\em Journal of Symbolic Computation}, 17:487--511, 1994.

\bibitem[WSV00]{WSV:00}
H.~Wolkowicz, R.~Saigal, and L.~Vandenberghe, editors.
\newblock {\em Handbook of Semidefinite Programming, Theory, Algorithms, and
  Applications}.
\newblock Kluwer Academic Publishers, 2000.

\bibitem[XB04]{XiB:04}
L.~Xiao and S.~Boyd.
\newblock Fast linear iterations for distributed averaging.
\newblock {\em Systems and Control Letters}, 53:65--78, 2004.

\bibitem[XBK07]{XBK:07}
L.~Xiao, S.~Boyd, and S.-J. Kim.
\newblock Distributed average consensus with least-mean-square deviation.
\newblock {\em Journal of Parallel and Distributed Computing}, 67:33--46, 2007.

\bibitem[YFK03]{SDPA6}
M.~Yamashita, K.~Fujisawa, and M.~Kojima.
\newblock Implementation and evaluation of {SDPA}~6.0 (semidefinite programming
  algorithm~6.0).
\newblock {\em Optimization Methods and Software}, 18:491--505, 2003.

\end{thebibliography}

\end{document}